\documentclass[12pt]{amsproc}
  \usepackage{latexsym} 
  \usepackage[all]{xy}
  \usepackage{amsfonts} 
  \usepackage{amsthm} 
  \usepackage{amsmath} 
  \usepackage{amssymb}
  \usepackage{pifont}

  \def\sw#1{{\sb{(#1)}}} 
   
  \def\su#1{{\sp{[#1]}}}

  \def\tens{\mathop{\otimes}}

  \def\<{{\langle}} 
  \def\>{{\rangle}} 
  \def\ra{{\triangleleft}} 
  \def\la{{\triangleright}} 
  \def\eps{\varepsilon}

  \def\note#1{{}} 
 \def\can{{\rm \textsf{can}}}

  \def\note#1{} 
  \def\M{{\bf M}} 
  \def\N{\mathbb{N}}

  \def\cD{{\mathcal D}}

  \def\cO{{\mathcal O}}

  \def\lrhom#1#2#3#4{{{}\sb{#1}{\rm Hom}\sb{#2}(#3,#4)}} 
   
  \def\rhom#1#2#3{{{\rm Hom}\sb{#1}(#2,#3)}}

  \def\rend#1#2{{{\rm End}\sb{#1}(#2)}} 
   \def\raut#1#2{{{\rm Aut}\sb{#1}(#2)}}

  \def\Rrhom#1#2#3#4{{{\rm Hom}\sp{#1}\sb{#2}(#3,#4)}}

  \def\C{\mathbb{C}}
  \def\R{\mathbb{R}} 
  \def\Z{\mathbb{Z}}

  \def\beq{\begin{equation}} 
  \def\eeq{\end{equation}}

  \def\id{\mathrm{id}} 
  \def\im{{\rm Im}} 
   
  \def\fg{{\mathfrak g}}

  \def\ot{{\otimes}}

  \def\roM{\varrho^{M}}

  \def\hsi{\hat{\sigma}}
   \def\bsi{\bar{\sigma}}

  \def\roA{{\varrho^A}}

  \def\coker{\mathrm{coker}}

    \def\oan#1{\Omega^{#1} (A)}
    \def\oau{\Omega A}
    \def\oaun#1{\Omega^{#1}\!A}
    \def\oa{\Omega (A)}
    \def\obn#1{\Omega^{#1} (B)}
    \def\ob{\Omega (B)}
    \def\obun#1{\Omega^{#1}\!B}
    
\def\1{\mathbb{I}}
\def\oomega{\overline{\omega}}
\def\cX{\mathfrak{X}}

\def\Rp{\Gamma_{+}}
\def\Rm{\Gamma_{-}}
\def \sfb{\mathsf{b}}
\def \sfa{\mathsf{a}}

\def \sfw{\mathsf{w}}

\def \sfu{\mathsf{u}}
\def\k{\Bbbk}

  \newcounter{zlist} 
  \newenvironment{zlist}{\begin{list}{(\arabic{zlist})}{ 
  \usecounter{zlist}\leftmargin2.5em\labelwidth2em\labelsep0.5em 
  \topsep0.6ex
  \parsep0.3ex plus0.2ex minus0.1ex}}{\end{list}}

  \newcounter{blist}

  \newcounter{rlist}

\newcommand{\osigma}{\overline{\sigma}}

\def\stac#1{\raise-.2cm\hbox{$\stackrel{\displaystyle\otimes}{\scriptscriptstyle{#1}}$}}

\def\cten#1{\raise-.2cm\hbox{$\stackrel{\displaystyle\widehat{\otimes}}
{\scriptscriptstyle{#1}}$}}

\def\tr{\mathrm{Tr}}

  \def\Label#1{\label{#1}\ifmmode\llap{[#1] }\else 
  \marginpar{\smash{\hbox{\tiny [#1]}}}\fi} 
  \def\Label{\label}

  \newtheorem{proposition}{Proposition}[section]
  \newtheorem{lemma}[proposition]{Lemma} 
  \newtheorem{corollary}[proposition]{Corollary} 
  \newtheorem{theorem}[proposition]{Theorem} 

  \theoremstyle{definition} 
  \newtheorem{definition}[proposition]{Definition}
    
  \newtheorem{example}[proposition]{Example} 

  \theoremstyle{remark}

\begin{document} 

 \title{Non-commutative integral forms and twisted multi-derivations} 
 \author{Tomasz Brzezi\'nski}
 \address{ Department of Mathematics, Swansea University, 
  Singleton Park, \newline\indent  Swansea SA2 8PP, U.K.} 
  \email{T.Brzezinski@swansea.ac.uk}   
  \author{Laiachi El Kaoutit}
  \address{Departamento de \'Algebra, Facultad de Educaci\'on y Humanidades de Ceuta, \newline\indent Universidad de Granada, El Greco No.\ 10,
E-51002 Ceuta, Spain}
  \email{kaoutit@ugr.es}
 \author{Christian Lomp}
 \address{Departamento de Matem\'atica Pura, Universidade do Porto, Porto, Portugal}
 \email{clomp@fc.up.pt}
    \date{\today} 
  \begin{abstract} 
Non-commutative connections of the second type or {\em hom-connections} and associated integral forms are studied as generalisations of {\em right connections} of Manin. First, it is proven that the existence of hom-connections with respect to the universal differential graded algebra is tantamount to the injectivity, and that every finitely cogenerated injective module admits a hom-connection with respect to any differential graded algebra. The bulk of the paper is devoted to describing a method of constructing hom-connections from \emph{twisted multi-derivations}. The notion of a {\em free} twisted multi-derivation is introduced and the induced first order differential calculus is described. It is shown that any free twisted multi-derivation on an algebra $A$ induces a unique hom-connection on $A$ (with respect to the induced differential calculus $\oan 1$) that vanishes on the dual basis of $\oan 1$. To any flat hom-connection $\nabla$ on $A$ one associates a chain complex, termed a {\em complex of integral forms} on $A$.  The canonical cokernel morphism to the zeroth homology space is called a {\em $\nabla$-integral}. Examples of free twisted multi-derivations, hom-connections  and corresponding integral forms are provided by covariant calculi on Hopf algebras  (quantum groups).  The example of a flat hom-connection within the 3D left-covariant differential calculus on the quantum group $\cO_q(SL(2))$ is described in full detail. A descent of hom-connections to the base algebra of a faithfully flat Hopf-Galois extension or a principal comodule algebra is studied. As an example, a hom-connection on the standard quantum Podle\'s sphere $\cO_q(S^2)$ is presented.
In both cases the complex of integral forms is shown to be isomorphic to the de Rham complex, and the $\nabla$-integrals coincide with Hopf-theoretic integrals or invariant (Haar) measures. 
 \end{abstract} 
  \maketitle

\section{Introduction}
The notion of a {\em right connection} appears in \cite[Chapter~4]{Man:gau} as a means for introducing integral forms and defining the  {Berezin integral} on a supermanifold. The prompt for this new type of connection comes from an approach to supersymmetric $\cD$-modules of Penkov \cite{Pen:mod}. A  formulation of the theory of  right connections or {\em co-connections} in terms of differential operators between modules over a   commutative algebra is described in \cite{Vin:coc}. It is also argued there that right connections or co-connections are dual notions (in a suitable sense) to connections in differential geometry. 

Motivated by the adjoint relationship between tensor products and homomorphisms, the notion of a {\em hom-connection} for any differential graded algebra was introduced in \cite{Brz:con}. Connections in non-commutative geometry, studied at least from its inception in \cite{Con:non}, are maps from a module $M$ to $M$ tensored with one-forms that satisfy the Leibniz rule (thus a non-commutative connection is a generalisation of the covariant derivation). Hom-connections are maps with a domain in homomorphisms from one-forms to $M$ and with $M$ as a codomain, that are again required to satisfy (a suitable version of) the Leibniz rule.  Under suitable commutativity and finiteness assumptions hom-connections reduce to right connections or co-connections in classical (super)geometry.
Similarly to connections, to a {\em flat} hom-connection one can associate a chain complex. This is the complex of module-valued {\em integral forms} in non-commutative geometry. 

Any finite-dimensional compact oriented (super)manifold admits a flat right connection with certain uniqueness property; see \cite[Chapter~4\S 5]{Man:gau}, \cite{Vin:coc}. This right connection is  dual to the de Rham differential and can be used to establish an equivalence between categories of left and right connections and an isomorphism between de Rham and integral forms. More directly, if $D$ is the dimension of the manifold, then the corresponding right connection arises from the de Rham differential from $D-1$-forms to $D$-forms. 
The finite-dimensionality plays in this construction the most crucial role.

The aim of the present paper is to extend the aforementioned classical construction of right connections to non-commutative algebras
 with a view on singling out algebras corresponding to finite-dimensional non-commutative manifolds through the existence of (flat) hom-connections 
 together with an isomorphism between de Rham and integral forms. This can be understood as a contribution to a programme aimed at understanding what algebras   describe finite-dimensional spaces in non-commutative geometry, a programme which is recently gathering momentum; see for example \cite{Kra:Poi}, \cite{Dub:non}.

The paper is organised as follows. In Section~\ref{sec.con.inj} we introduce the notions of integral forms  and integrals associated to a hom-connection, and consider existence of hom-connections. In particular, it  is  proven that  existence of hom-connections with respect to the universal differential graded algebra is tantamount to the injectivity, and that every finitely cogenerated injective module admits a hom-connection with respect to any differential graded algebra. Section~\ref{sec.twist.deriv} describes a construction of hom-connections from free twisted multi-derivations. First the notion of a {\em free} twisted multi-derivation is introduced and the induced first order differential calculus is described. Sufficient and necessary conditions for a twisted multi-derivation with a triangular twisting matrix to be free are determined in Proposition~\ref{triangular}.  It is shown that any free multi-derivation on an algebra $A$ induces a unique hom-connection on $A$ (w
 ith respect to the induced differential calculus $\oan 1$) that vanishes on 
 the dual basis of $\oan 1$; see Theorem~\ref{thm.hom-der} and Corollary~\ref{cor.main}. Examples in this section  include hom-connections on quasi-free or smooth algebras and on algebras with a derivation based differential calculus. In particular, the $\nabla$-integral on a  matrix algebra $M_n(\C)$ with derivation based differential calculus is shown to coincide with the integral constructed in \cite{DubKer:non}. As another very explicit illustration of Theorem~\ref{thm.hom-der} we construct a hom-connection and complex of integral forms for the quantum plane with a covariant differential structure. We show that this complex is isomorphic to the non-commutative de Rham complex.
 
 Main examples are described in Section~\ref{sec.examples}. These are examples of hom-connections and integral forms on algebras with a Hopf-algebra coaction including Hopf algebras themselves. It is proven that any left- (or right-) covariant  differential calculus on a Hopf algebra with bijective antipode gives rise to a free twisted multi-derivation and hence to a hom-connection; see Theorem~\ref{thm.cov}. Furthermore a close relationship between integrals on Hopf algebras and integrals associated to this hom-connection is established. A hom-connection within the Woronowicz 3D left-covariant differential calculus on the quantum group $\cO_q(SL(2))$ is described in full detail. Next the descent of hom-connections within Hopf-Galois extensions is studied. It is shown that a hom-connection on a total algebra $A$ of a faithfully flat Hopf-Galois extension  $B\subseteq A$ descends to a
  hom-connection on $B$, provided it has a certain covariance property and the (strongly) horizontal forms on $B\subseteq A$ are a direct summand (in a suitable category) of forms on $A$; see Theorem~\ref{thm.induction}. In particular, a covariant hom-connection on a principal comodule algebra $B\subseteq A$ with respect to the universal differential calculus descends to a hom-connection on $B$. Also, a hom-connection on $\cO_q(SL(2))$ (with respect to the 3D calculus) described earlier descends to a hom-connection on the quantum standard Podle\'s sphere $\cO_q(S^2)$. In both cases,  hom-connections can be identified with  exterior differentials $\Omega^2(\cO_q(SL(2))) \to \Omega^3(\cO_q(SL(2)))$ and $\Omega^1(\cO_q(S^2)) \to \Omega^2(\cO_q(S^2))$, respectively. Similarly to the case of the quantum plane, the constructed hom-connections are flat and the corresponding complexes of integral forms are isomorphic to de Rham complexes $\Omega(\cO_q(SL(2)))$ and $\Omega(\cO_q(S^2))$
 , respectively.  This is a non-commutative counterpart of the classical identification of differential and integral forms on a compact oriented finite dimensional manifold.  

We work over a field $\k$. All algebras are unital and associative. Only Section~\ref{sec.examples} requires the reader to have some familiarity with the language of Hopf algebras or quantum groups. We follow the standard Hopf algebra notation and conventions there.

\section{Hom-connections and injective modules}\label{sec.con.inj}
\subsection{Hom-connections, integral forms and gauge transformations} By a {\em differential graded algebra over an algebra $A$} we mean a non-negatively graded differential graded algebra $(\oa, d)$ with $\oan 0 = A$. The pair $(\oan 1, d)$ is referred to as a {\em first order differential calculus} on $A$. Any algebra $A$ admits the {\em universal differential graded algebra $(\oau, d)$} over $A$ defined as follows. $\oau: = T_A(\oaun 1)$ is the tensor algebra of the $A$-bimodule $\oaun 1 = \ker \mu$, where $\mu: A\ot A\to A$ is the multiplication map. The differential is defined as $d: A\to \oaun 1$, $a\mapsto 1\ot a - a\ot 1$ and extended to the whole of $\oau$ by the graded Leibniz rule.

A right {\em hom-connection} with respect to  a differential graded algebra $(\oa, d)$ over an algebra $A$, is a pair $(M,\nabla)$, where $M$ is a right $A$-module and 
$$
\nabla :\rhom A {\oan 1} M \to M,
$$
is a $\k$-linear map, such that, for all $f\in \rhom A {\oan 1} M$ and $a\in A$,
$$
\nabla (fa) = \nabla (f)a + f(da).
$$
Here $\rhom A {\oan 1} M $ is a right $A$-module by $(fa)(\omega) := f(a\omega)$, $\omega \in \oan 1$.

Any hom-connection $(M,\nabla )$ can be extended to higher forms. The vector space $\bigoplus_{n=0} \rhom A {\oan n} M$ is a right module over $\oa$ with the multiplication, for all $\omega\in \oan n$, $f\in \rhom A{\oan{n+m}} M$, $\omega'\in \oan m$,
$$
f\omega (\omega') := f(\omega\omega').
$$
For any $n>0$, define 
$
\nabla_n: \rhom A {\oan {n+1}}M \to \rhom A {\oan {n}}M,
$
by
\begin{equation}\label{nablan}
\nabla_n(f)(\omega) :=  \nabla  (f\omega) + (-1)^{n+1} f(d\omega),
\end{equation}
for all $f\in \rhom A {\oan{n+1}}M$ and $\omega \in \oan{n}$. 

The map $F := \nabla \circ \nabla_1$ is called the {\em curvature} of  $(M,\nabla )$, and  $(M,\nabla )$ is said to be {\em flat} provided $F=0$. To a flat hom-connection $(M,\nabla )$ one associates a chain complex $\left(\bigoplus_{n=0} \rhom A {\oan n} M, \nabla\right)$. The homology of this complex is denoted by $H_*(A; M,\nabla)$. In case $M=A$ this complex is termed a {\em complex of integral forms} on $A$, and the canonical map $$\Lambda: A \longrightarrow  \coker (\nabla) = H_0(A; A,\nabla)$$ is called a {\em $\nabla$-integral on $A$}.

We complete this preliminary section by describing the action of the group of module automorphisms on hom-connections. This is a hom-connection version of gauge transformations of connections. 

\begin{proposition}\label{prop.gauge}
Let  $G= \raut A M$ be the group of automorphisms of a right $A$-module $M$, and view  $\rhom A {\oan 1} M$ as a right $G$-space by 
$$
f \ot \Phi \mapsto \left[f\ra\Phi: \omega \mapsto \Phi^{-1}\left(f\left(\omega\right)\right)\right].
$$
For any hom-connection $(M,\nabla )$ and $\Phi \in \raut A M$, the pair $(M,\nabla ^\Phi)$,
where
$$ 
\nabla ^\Phi : \rhom A {\oan 1} M\to M, \qquad  f \mapsto \Phi\left( \nabla \left( f\ra\Phi\right)\right),
$$
is a hom-connection.
\end{proposition}
\begin{proof}
Note that, for all $f\in \rhom A {\oan 1} M$, $\Phi\in \raut AM$ and $a\in A$,
$$
(fa)\ra \Phi = (f\ra \Phi)a.
$$
Therefore,
\begin{eqnarray*}
\nabla^\Phi(fa) &=& \Phi\left( \nabla \left( (fa)\ra \Phi\right)\right)  = \Phi\left( \nabla \left( (f\ra\Phi)a\right)\right)
= \Phi\left( \nabla \left( f\ra\Phi\right) a\right) + \Phi(f\ra\Phi(da)) \\
&=&  \nabla^\Phi(f)a + \Phi\left(\Phi^{-1}\left( f\left( da\right)\right)\right) = \nabla^\Phi(f)a +f(da),
\end{eqnarray*}
where the fourth equality follows by the right $A$-linearity of $\Phi$. Hence $\nabla^\Phi$ is a hom-connection as claimed.
\end{proof}

The action of $G= \raut A M$ extends to the action on $\bigoplus_{n=0} \rhom A {\oan n} M$. A similar calculation to the proof of Proposition~\ref{prop.gauge} yields the following relation between the curvatures of $\nabla$ and $\nabla^\Phi$, 
$$
F^\Phi(f) = \Phi\left(F\left(f\ra\Phi\right)\right), 
$$
for all $f\in \rhom A {\oan 2} M$. In particular, the gauge transform of a flat hom-connection is a flat hom-connection. 

\subsection{Existence of hom-connections}
First we determine, when a module admits a hom-connection with respect to the universal differential calculus, and thus obtain a dual version of \cite[Corollary~8.2]{CunQui:alg}.

\begin{theorem}\label{thm.inj}
A right $A$-module $M$ admits a hom-connection with respect to the universal differential graded algebra if and only if $M$ is an injective module.
\end{theorem}
\begin{proof}
Note that $M$ is an injective right $A$-module if and only if there is a right $A$-module map $\pi: \rhom \k A M \to M$ such that $\pi\circ \theta = {\rm id}$, where $\theta: M\to \rhom \k A M$ is the canonical monomorphism induced by the $A$-multiplication on $M$, i.e.\  $\theta(m)(a) = ma$, for all $m\in M$ and $a\in A$. Indeed, by the standard argument $\rhom \k A M$ is injective, hence if $\pi$ exists, $M$ is a direct summand of an injective module, therefore it is injective. Conversely, if $M$ is an injective right $A$-module, then the following diagram
$$
\xymatrix{ 0 \ar[r] & M\ar[r]^-\theta \ar[d]_= & \rhom \k  AM\\
& M, & }
$$
can be completed by the required right $A$-module map $\pi: \rhom \k AM \to M$. 

In the light of \cite[Section~3.9]{Brz:con}, hom-connections on $M$ with respect to the universal differential graded algebra are in bijective correspondence with right $A$-module maps $\phi:  \rhom A {A\ot A} M \to M$ such that
\begin{equation}\label{eq.phi}
\phi\circ \rhom A \mu M \circ \psi = {\rm id},
\end{equation}
where $\psi: M\to \rhom A AM$ is the canonical isomorphism $m\mapsto [a\mapsto ma]$. Write $\bar{\psi} : \rhom \k AM \to \rhom A {A\ot A} M$ for the canonical isomorphism $f\mapsto [a\ot b \mapsto f(a)b]$, and note that
$$
\theta = \bar{\psi}^{-1}\circ \rhom A \mu M\circ \psi.
$$
If a hom-connection, i.e.\ $\phi$ satisfying \eqref{eq.phi} exists, then define
$
\pi = \phi\circ \bar{\psi},
$ 
and compute
$$
\pi\circ \theta = \phi\circ \bar{\psi}\circ \bar{\psi}^{-1}\circ \rhom A \mu M\circ \psi  = {\rm id}.
$$
 Therefore, $M$ is an injective right $A$-module. 

Conversely, assume that $M$ is an injective right $A$-module with the corresponding $\pi$, and define
$
\phi = \pi \circ \bar{\psi}^{-1}.
$
Then 
$$
\phi\circ  \rhom A \mu M \circ \psi = \pi \circ \bar{\psi}^{-1} \circ  \rhom A \mu M \circ \psi = \pi\circ\theta = {\rm id},
$$
hence there is a hom-connection in $M$.
\end{proof}

With no restriction on the differential graded algebra, 
\begin{proposition}\label{prop.dir.sum}
A finite direct sum of $A$-modules admits a hom-connection if and only if each of the summands admits a hom-connection.
\end{proposition}
\begin{proof}
Suppose that $M= N\oplus P$ and let  $\pi_N: M\to N$ and $\iota_N: N\to M$ be the (direct-sum defining) epimorphism and monomorphism. If $M$ has a hom-connection $\nabla ^M :\rhom A {\oan 1} M\to M$, then 
$$
\nabla ^N :\rhom A {\oan 1} N\to N, \qquad f\mapsto \pi_N\left(\nabla ^M\left(\iota_N\circ f\right)\right),
$$
is a hom-connection on $N$.

Conversely, assume that $(N,\nabla ^N)$ and $(P,\nabla ^P)$ are hom-connections. Let $\pi_P: M\to P$ and $\iota_P: P\to M$ be the (direct-sum defining) epimorphism and monomorphism. Then the hom-connection 
$$
\nabla ^M: \rhom A {\oan 1} {N\oplus P} \simeq \rhom A {\oan 1} N \oplus \rhom A {\oan 1} {P} \to N\oplus P,
$$
is given by
$$
\nabla ^M(f) = \iota_N\left(\nabla^N\left(\pi_N\circ f\right)\right) + \iota_P\left(\nabla^P\left(\pi_P\circ f\right)\right).
$$
In both cases checking that the maps given by the stated formulae satisfy the Leibniz rule for a hom-connection is a routine calculation.
\end{proof}

The following proposition is a dual version of \cite[Proposition~III.3.6]{Con:ncg}.
\begin{proposition}\label{prop.cogen}
Every finitely cogenerated injective right $A$-module admits a hom-connection (with respect to any differential graded algebra over $A$).
\end{proposition}
\begin{proof}
Recall from \cite{Vam:dua} that an injective right $A$-module $M$ is finitely cogenerated if and only if it is a finite direct sum of injective envelopes of simple $A$-modules, i.e.\ $M = \oplus_{i=1}^n E(S_i)$. Since every $S_i$ is a simple module, $E(S_i)$ is a direct summand of $A^* = \rhom \k  A\k $, where $A^*$ is a right $A$-module by the left multiplication of arguments, i.e.\ $(\xi a)(b) = \xi(ab)$. $A^*$ has a hom-connection $-d^* := - \rhom \k  d \k : \rhom A{\oan 1} {A^*} \to A^*$; see  \cite[Section~2.8]{Brz:con}. Hence, by Proposition~\ref{prop.dir.sum}, all $E(S_i)$, $i=1,2,\ldots, n$ have hom-connections, and so does their direct sum $M$.
\end{proof}

\section{Twisted multi-derivations and hom-connections}\label{sec.twist.deriv}
\setcounter{equation}{0}
By a {\em right twisted multi-derivation} in an algebra $A$ we mean a pair $(\partial, \sigma)$, where $\sigma: A\to M_n(A)$ is an algebra homomorphism ($M_n(A)$ is the algebra of $n\times n$ matrices with entries from $A$) and $\partial :A\to A^n$ is a $\k $-linear map such that, for all $a\in A$, $b\in B$,
\begin{equation}\label{eq.partial}
\partial(ab) = \partial(a)\sigma(b) + a\partial(b).
\end{equation}
Here $A^n$ is understood as an $(A$-$M_n(A))$-bimodule.
Thus, writing $\sigma(a) = (\sigma_{ij}(a))_{i,j=1}^n$ and $\partial(a) = (\partial_i(a))_{i=1}^n$, 
 \eqref{eq.partial} is equivalent to the following $n$ equations
$$
\partial_i(ab) = \sum_j \partial_j(a)\sigma_{ji}(b) + a\partial_i(b), \qquad i=1,2,\ldots, n.
$$
Equivalently, a right twisted multi-derivation is a derivation on $A$ with values in the $A$-bimodule 
$A^n_\sigma$ with the underlying vector space equal to $A^n$, but with $A$-multiplications
\begin{equation}\label{def.sigma.a}
a \cdot x \cdot c =  (a\sum_jx_j\sigma_{ji}(c)), \qquad \forall a,c\in A,\; x=(x_i)\in A^n_\sigma.
\end{equation}
Twisted multi-derivations (for $n=2$) give rise to {\em double Ore extensions} \cite{ZhaZha:dou}.

A map $\sigma : A \to M_n(A)$ can be  viewed as an element of $M_n(\rend \k  A)$. Write $\bullet$ for the product in $M_n(\rend \k  A)$,  $\1$ for the unit in $M_n(\rend \k  A)$ and $\sigma^T$ for the transpose of $\sigma$. 
\begin{definition}\label{def.der.free}
Let $(\partial, \sigma)$ be a right twisted multi-derivation. We say that $(\partial, \sigma)$ is {\em free}, provided there exist algebra maps $\bsi: A \to M_n(A)$ and $\hsi: A \to M_n(A)$ such that
\begin{equation}\label{bar.sigma}
\bsi \bullet \sigma^T = \1, \qquad \sigma^T\bullet\bsi = \1\, ,
\end{equation}
\begin{equation}\label{hat.sigma}
\hsi \bullet \bsi^T = \1, \qquad \bsi^T\bullet\hsi = \1\, .
\end{equation}
\end{definition}

For any algebra homomorphism $\sigma: A \rightarrow M_n(A)$ define  the $A$-bimodule $A^n_\sigma$ as free left $A$-module $\bigoplus_{i=1}^n A \omega_i$  with basis $\omega_1,\ldots, \omega_n$ and right $A$-action given by 
\begin{equation}\label{eq.dif.rel}
 \omega_i a = \sum _j \sigma_{ij}(a) \omega_j ,  \qquad i=1,2,\ldots, n.
\end{equation}
Identifying an arbitrary element $x\in A^n_\sigma$ with the coefficient vector $(x_i)$ of its representation using the left $A$-basis $\omega_1,\ldots, \omega_n$ we obtain the formula \eqref{def.sigma.a}.

\begin{lemma}
There exists an algebra homomorphism $\osigma:A\rightarrow M_n(A)$ such that 
$$\sigma^T \bullet \osigma = \mathbb{I} = \osigma \bullet \sigma^T,$$
if and only if $\omega_1, \ldots, \omega_n$ is also a basis for $A^n_\sigma$ as a right $A$-module.
\end{lemma}
\begin{proof}
Suppose that $\osigma$ exists, then $\sum_k \osigma_{jk}(\sigma_{ik}(a)) = \delta_{ji} a$, for all $a\in A$.
Suppose that $ \sum_i \omega_i a_i = 0$, for some $a_i\in A$. Then 
$ \sum_{ij} \sigma_{ij}(a_i) \omega_j = 0$. Since $\{ \omega_i \}_i$ is  a free basis of the left $A$-module, 
$$ 
\sum_i\sigma_{ij}(a_i) = 0, \qquad \mbox{for all } j=1, 2,\dots , n.
$$
In particular $\sigma_{kj}(a_k) = - \sum_{i\neq k} \sigma_{ij}(a_i)$ for all $j$. Hence, for all $k$,
$$ a_k = \sum_{j=1}^n \osigma_{kj}\circ \sigma_{kj}(a_k) = - \sum_{j=1}^n\sum_{i\neq k} \osigma_{kj}\circ\sigma_{ij}(a_i) = -\sum_{i\neq k} \delta_{ik} a_i = 0.$$
Therefore, $\sum_i \omega_iA = \bigoplus \omega_iA$ is free on the right. As 
$$ \sum_j \omega_i\osigma_{ij}(a) = \sum_{j,k} \sigma_{ik}(\osigma_{ij}(a))\omega_k = a\omega_i,$$
it follows  that $A^n_\sigma=\bigoplus \omega_iA$.

On the other hand suppose that $\{ \omega_i\}_i$ is a right $A$-basis for $A_\sigma^n$. For all $i= 1,2,\ldots, n$ and $a\in A$,  $a\omega_i = \sum_j\omega_j\osigma_{ji}(a)$, for some well-defined $\osigma_{ji}(a)\in A$. Since $A^n_\sigma$ is a left $A$-module and free on the right  with the basis $\{ \omega_i\}_i$, the function $\osigma: a \mapsto (\osigma_{ij}(a)) \in M_n(A)$ is an algebra map.
Then $a\omega_i = \sum_ j\omega_j(\osigma_{ji}(a)) = \sum_{j,k}\sigma_{jk}\osigma_{ji}(a)\omega_k$, and the freeness on the left shows that 
$\delta_{ki}a = \sum_j\sigma_{jk}(\osigma_{ji}(a))$, i.e.\ $\sigma^T \bullet \osigma = \mathbb{I}$. The other equality follows analogously.
\end{proof}

In other words a right twisted multi-derivation $(\partial, \sigma)$ is free if and only if $\omega_1, \ldots, \omega_n$ is a right $A$-basis for $A_\sigma^n$ and $A_{\osigma}^n$.

In general it might be difficult to determine, when an algebra map $\sigma : A\to M_n(A)$ admits maps $\bsi$, $\hsi$ satisfying conditions of Definition~\ref{def.der.free} (note that we are dealing with matrices with non-commutative entries even if $A$ is a commutative algebra). However, when $\sigma(a)$ is an upper-triangular matrix for each $a \in A$, a criterion for existence of $\bsi$, $\hsi$ can be established.

\begin{proposition}\label{triangular} 
Let $(\partial,\sigma)$ be a right twisted multi-derivation in which  $\sigma$ is an upper-triangular matrix (in $M_n(\rend \k  A)$) with non zero diagonal entries.  Then, $(\partial, \sigma)$ is a free multi-derivation in the sense of Definition \ref{def.der.free} with $\bsi$ a lower-triangular matrix if and only if the  diagonal entries $\sigma_{i\,i}$, $i=1,\cdots, n$ of $\sigma$ are invertible.
\end{proposition}
\begin{proof}
To ease the notation we denote the composition of endomorphisms by juxtaposition. If $(\partial, \sigma)$ is a free multi-derivation with $\bsi$ a lower-triangular matrix, then equations \eqref{bar.sigma} immediately imply that, for all $i$,   $\sum_k \sigma_{ki}\bsi_{ki} = \sigma_{ii} \bsi_{ii} =1$, and  $\sum_k \bsi_{ik}\sigma_{ik} = \bsi_{ii} \sigma_{ii} = 1$. Thus all diagonal entries in $\sigma$ are invertible. \\ Conversely assume that all diagonal entries in $\sigma$ are invertible and define a  matrix $\bsi$ by
\begin{eqnarray*}
 \bsi_{i\,j} & =& 0, \quad \text{for } i<j; \\
\bsi_{i\,i} & =& \sigma_{i\,i}^{-1}, \quad \text{for all } i=1,\cdots,n;\\
\bsi_{i\,j} & =& -\sum^{i-1}_{k=j}\sigma_{i\,i}^{-1} \sigma_{k\, i} \bsi_{k\, j} , \quad \text{for } i\geq j+1.
\end{eqnarray*}
We claim that $\bsi \bullet \sigma^T =  \sigma^T\bullet\bsi = \1\,$. Clearly, $\sum_k \sigma_{ki}\bsi_{ki} = \sigma_{ii} \bsi_{ii} =1$, and similarly $\sum_k \bsi_{ik}\sigma_{ik} =1$. It remains to prove that $\sum_k \sigma_{k\, j}\bsi_{k\, i} = \sum_k \bsi_{i\, k}\sigma_{j\,k}=0$, for any $i\neq j$. Given such a pair of indices, it is clear from definitions that $\sum_k \sigma_{k\, j}\bsi_{k\, i}=0$, for any $j < i$. On the other hand, if $j > i$, then 
\begin{eqnarray*}
\sum_{k=1}^n \sigma_{k\, j} \bsi_{k\, i} &=&  \sum_{k=i}^j \sigma_{k\, j} \bsi_{k\, i}  \,=\,   \sum_{k=i}^{j-1} \sigma_{k\, j} \bsi_{k\, i} + \sigma_{j\, j} \bsi_{j\, i}\\
& = &  \sum_{k=i}^{j-1} \sigma_{k\, j} \bsi_{k\, i} - \sigma_{i\, j} \sigma^{-1}_{i\, i} - \sum_{k=i+1}^{j-1} \sigma_{k\, j}\bsi_{k\, j} \,=\, 0. 
\end{eqnarray*}
For the second equality in \eqref{bar.sigma}, observe that $\sum_k\bsi_{i\, k}\sigma_{j\, k}=0$, for $i<j$. So let $i > j$. If $i=j+1$, then,
  by the definition of $\bsi$, 
  $$ 
  \sum_{k=1}^n \bsi_{i\, k}\sigma_{j\, k}\,=\, \sum_{k=j}^{j+1}\bsi_{i\, k}\sigma_{j\, k}\,=\, \bsi_{j+1\,j}\sigma_{j\,j} + \bsi_{j+1\, j+1}\sigma_{j\,j+1} \,=\,0.
  $$ 
  Assume that $\sum_{k=j}^l\bsi_{l\, k}\sigma_{j\, k}\,=\,0$, for all  $l$ such that $j< l\leq i$ and $l-j \geq 2$. Then 
\begin{eqnarray*}
\sum_{k=1}^n \bsi_{i+1\, k}\sigma_{j\, k} &=& \sum_{k=j}^{i+1} \bsi_{i+1\, k}\sigma_{j\, k}  
= \sum_{k=j}^{i} \bsi_{i+1\, k}\sigma_{j\, k} + \bsi_{i+1\, i+1}\sigma_{j\, i+1} \\ &=& -\sum_{k=j}^i \sum_{l=k}^i \sigma_{i+1\,i+1}^{-1}\sigma_{l\, i+1} \bsi_{l\, k}\sigma_{j\,k} +\bsi_{i+1\, i+1}\sigma_{j\, i+1} \\ &=& -\sum_{l=j}^i \sigma_{i+1\,i+1}^{-1}\sigma_{l\, i+1}\left(\underset{}{} \sum_{k=j}^l \bsi_{l\, k}\sigma_{j\,k}\right) +\bsi_{i+1\, i+1}\sigma_{j\, i+1} \\ 
&=& -\sigma_{i+1\,i+1}^{-1}\sigma_{j\,i+1}(\bsi_{j\, j}\sigma_{j\,j})  -\left(\underset{}{} \sum_{l=j+1}^i \sigma_{i+1\,i+1}^{-1}\sigma_{l\, i+1}\left(\underset{}{} \sum_{k=j}^l \bsi_{l\, k}\sigma_{j\,k}\right) \right) \\
&& + \bsi_{i+1\, i+1}\sigma_{j\, i+1} \\ 
&=& -\sigma_{i+1\,i+1}^{-1}\sigma_{j\,i+1}(\sigma_{j\, j}^{-1}\sigma_{j\,j})+\sigma_{i+1\, i+1}^{-1}\sigma_{j\, i+1} \,\,=\,\, 0\, ,
\end{eqnarray*}
where in the third equality  the definition $\bsi_{i+1\, k}= -\sum_{l=k}^i \sigma_{i+1\,i+1}^{-1}\sigma_{l\,i+1}\bsi_{l\, k}$ and in the fifth one the induction hypothesis were used. This finishes the proof of equations \eqref{bar.sigma}.
In a similar way, one can show that if $\alpha \in M_n({\rm End}_{\k }(A))$ is a lower-triangular matrix in which all the $\alpha_{i\,i}$ are invertible, then the upper-triangular matrix $\widetilde{\alpha}$ defined by 
$$
\widetilde{\alpha}_{i\,i}=\alpha_{i\,i}^{-1},\text{ for } i=1,\cdots,n,\,\text{ and } \widetilde{\alpha}_{i\,j}=-\sum_{l=i+1}^j\alpha_{i\,i}^{-1}\alpha_{l\,i}\widetilde{\alpha}_{l\,j}, \text{ for } i+1 \leq j,
$$ 
satisfies $\widetilde{\alpha} \bullet \alpha^T \,=\, \alpha^T \bullet \widetilde{\alpha} \,=\, \1$. In particular, if  $\hat{\sigma} :=\widetilde{\bsi}$, then $ \hat{\sigma} \bullet \bsi^T\, =\, \bsi^T \bullet \hat{\sigma}\,=\, \1$, i.e.\ equations \eqref{hat.sigma} are satisfied as required. 
\end{proof}

Given a free right twisted multi-derivation $(\partial, \sigma; \bsi, \hsi)$  define a first order differential calculus on $A$ as follows: 
$
\oan 1 =  A^n_{\sigma},
$
hence writing $\omega_i = (0,0,\ldots, 1,\ldots ,0)$ for standard generators of $A^n_{\sigma}$ we obtain relations  \eqref{eq.dif.rel} in $\oan 1$.
Equations \eqref{eq.dif.rel} and the freeness of $\oan 1$ as an $A$-module imply that
\begin{equation}\label{eq.a.om}
a\omega_i = \sum_j \omega_j\bsi_{ji}(a),  \qquad i=1,2,\ldots, n.
\end{equation}
The exterior differential $d:A \to \oan 1$ is defined by
\begin{equation}\label{eq.dif.d}
da = \sum_i \partial_i(a)\omega_i = \sum_{i,j}\omega_i \bsi_{ij}\left(\partial_j\left(a\right)\right).
\end{equation}
The fact that  $(\partial, \sigma)$ is a right twisted multi-derivation ensures that the differential $d$ satisfies the Leibniz rule. This first order differential calculus can be extended to a differential graded algebra in a standard way. Furthermore, $\oan 1$ can be understood as a universal calculus in which relations \eqref{eq.dif.rel}--\eqref{eq.dif.d} are satisfied.

Within this set-up, one can formulate the following non-commutative version of \cite[Chapter~4\S 5, Proposition~3]{Man:gau}.
\begin{theorem}\label{thm.hom-der}
Let $(\partial, \sigma; \bsi ,\hsi)$,  be a free right twisted multi-derivation on $A$, and let $\oan 1$ be the associated first order differential calculus with generators $\omega_i$. Define right $A$-module maps $\xi_i : \oan 1\to A$ by $\xi_i(\omega_j) = \delta_{ij}$, $i,j = 1,2,\ldots , n$. Then there exists a unique hom-connection $\nabla : \rhom A {\oan 1} A\to A$ such that $\nabla (\xi _i) =0$, for all $i = 1,2, \ldots , n$. 
\end{theorem}
\begin{proof}
For each $i=1,2,\ldots, n$, write $\partial_i^\sigma := \sum_{j,\,k} \bsi_{kj}\circ\partial_j\circ \hsi_{ki}$, and define
\begin{equation}\label{def.hom.twist}
\nabla : \rhom A {\oan 1} A\to A, \qquad f\mapsto \sum_i \partial_i^\sigma \left( f\left(\omega_i\right)\right)\, .
\end{equation}
Since $\partial_i^\sigma(1) =0$, $\nabla (\xi _i) =0$, for all $i = 1,2, \ldots , n$. Furthermore, for all $a\in A$ and $f\in \rhom A {\oan 1} A$,
\begin{eqnarray*}
\nabla (fa) &=& \sum_ i \partial_i^\sigma\left( f\left(a\omega_i\right)\right) = \sum_{i,\, j,\, k, \,l}\bsi_{kj}\left( \partial_j\left(\hsi_{ki}\left( f\left(\omega_l\right) \bsi_{li}\left(a\right)\right)\right)\right)\\
&=&  
\sum_{j,\, k,\,l}\bsi_{kj}\left( \partial_j\left(\hsi_{kl}\left( f\left(\omega_l\right)\right) a\right)\right) \\
&=& \sum_{j,\,k,\,l,\,r} \bsi_{kj}\left( \partial_r\left(\hsi_{kl}\left( f\left(\omega_l\right)\right)\sigma_{rj}\left( a\right)\right)\right)
+ \sum_{j,\,k,\,l}\bsi_{kj}\left( \hsi_{kl}\left( f\left(\omega_l\right)\right) \partial_j\left(a\right)\right)\\
&=& \sum_i\partial_i^\sigma\left( f\left(\omega_i\right)\right)a +   \sum_{j,\,l}f\left(\omega_l\right)\bsi_{lj}\left(\partial_j\left( a\right)\right)\\
&=& \nabla (f)a + \sum_{j,\,l} f\left(\omega_l\bsi_{lj}\left(\partial_j\left( a\right)\right)\right) = \nabla (f)a + f(da),
\end{eqnarray*}
where the first  equality follows by the definition of right $A$-action on $\rhom A {\oan 1} A$ and the definition of $\nabla $ in \eqref{def.hom.twist}. The second equality follows by  relations \eqref{eq.a.om} and the $A$-linearity of $f$. The third equality is a consequence of multiplicativity of $\hsi$ and the first of equations \eqref{hat.sigma}. The fourth equality follows by the twisted derivation property, while the fifth  one is a consequence of multiplicativity of $\bsi$, the first of equations \eqref{bar.sigma} and the second of \eqref{hat.sigma}. Finally the linearity of $f$ and  the definition of $d$ in \eqref{eq.dif.d} are used. Hence $\nabla $ is a hom-connection. 

Similar calculations (that in particular use the second of equations \eqref{hat.sigma}) and the definition of $\xi_i$ affirm that, for all $f\in \rhom A {\oan 1} A$,
$$
f =  \sum_{i,\,k} \xi_i \hsi_{ik}\left( f\left(\omega_k\right)\right).
$$
Suppose that $\bar{\nabla}$ is a hom-connection such that $\bar{\nabla}(\xi_i) =0$, $i=1,2,\ldots ,n$. Then
\begin{eqnarray*}
\bar{\nabla}(f) &=& \sum_{i,\,k} \bar{\nabla}\left(\xi_i \hsi_{ik}\left( f\left(\omega_k\right)\right)\right) = \sum_{i,\,k} \bar{\nabla}\left(\xi_i \right)\hsi_{ik}\left( f\left(\omega_k\right)\right)  + \sum_{i,\,k} \xi_i \left( d\hsi_{ik}\left( f\left(\omega_k\right)\right)\right)\\
&=& \sum_{i,\,j,\,k,\,l} \xi_i \left( \omega_l \bsi_{lj}\left(\partial_j\left(\hsi_{ik}\left( f\left(\omega_k\right)\right)\right)\right)\right) = \sum_{i,\,j,\,k} \bsi_{ij}\left(\partial_j\left(\hsi_{ik}\left( f\left(\omega_k\right)\right)\right)\right) = \nabla (f).
\end{eqnarray*}
The second equality is the Leibniz rule for a hom-connection, the third one is a consequence of the hypothesis on $\bar{\nabla}$ and the definition of $d$ in \eqref{eq.dif.d}. The remaining two equalities are consequences of the definitions of the $\xi_i$ and $\nabla $. This completes the proof of uniqueness of $\nabla $.
\end{proof}

A first order differential calculus $\oan 1$ is said to be {\em dense} if every element of $\oan 1$ is of the form $\sum_i a_idb_i$, for some $a_i, b_i\in A$. The calculus discussed in Theorem~\ref{thm.hom-der} is dense if and only if there exist two finite subsets $\{a_{i\,t}\}, \{b_{i\,t}\}$ of elements of $A$  such that, 
$$
\sum_t a_{i\,t}\partial_k(b_{i\,t}) = \delta_{ik}, \quad \text{for all } i,k =1, \ldots, n.
$$
Typically, one is interested in calculi that are dense. Also typically the calculi of main interest are not free as $A$-modules (but they are often finitely generated and projective, as $A$ is understood as functions on a non-commutative space and $\oan 1$ is understood as sections of the non-commutative cotangent bundle). 
Note, however, that the calculations in the proof of Theorem~\ref{thm.hom-der} justify the following assertion. 

\begin{corollary}\label{cor.main}
Let $(\partial, \sigma; \bsi, \hsi)$ be a free right-twisted multi-derivation, and let $\oan 1$ be {\em any} differential calculus on $A$ finitely (but not necessarily freely) generated by the $\omega_i$ and such that relations \eqref{eq.dif.rel}--\eqref{eq.dif.d} are satisfied, for all $a\in A$.  Then the formula \eqref{def.hom.twist} defines a hom-connection on $A$ with respect to $\oan 1$. 
\end{corollary}

In a typical non-commutative geometry situation, rather than constructing calculus from a multi-derivation, one would start with a suitable differential calculus, and then search for a free right-twisted multi-derivation. In all such situations it is Corollary~\ref{cor.main} rather than Theorem~\ref{thm.hom-der} that produces a hom-connection.

The construction in Theorem~\ref{thm.hom-der} or Corollary~\ref{cor.main} simplifies if $(\delta, \sigma)$ is a right twisted multi-derivation with a diagonal matrix $\sigma$. Write $\sigma_i$ for the (only non-zero) diagonal elements of $\sigma$. Then each of the $\partial_i$ becomes a right twisted derivation, i.e. $\partial_i(ab) = \partial_i(a)\sigma_i(b) +a\partial_i(b)$. In this case the conditions \eqref{hat.sigma} and \eqref{bar.sigma} are mutually equivalent and simply state that each of the endomorphisms $\sigma_i$ is an algebra automorphism; see Proposition~\ref{triangular}. Furthermore $\bsi$ is the inverse of $\sigma$ (i.e.\ a diagonal matrix with entries $\sigma_i^{-1}$) and $\hsi = \sigma$. If there exist scalars $q_i$ such that $\sigma^{-1}_i\partial_i \sigma_i = q_i \partial_i $, then, following \cite{GooLet:pri}, each of the $\partial_i$ is called a {\em $q_i$-skew derivation}. In this case the formula for a hom-connection in Theorem~\ref{thm.hom-der} takes 
 particularly simple form: $$\nabla (f) = \sum_i q_i \partial_i (f(\omega_i)).$$ 
 
To illustrate the construction in Theorem~\ref{thm.hom-der} or Corollary~\ref{cor.main} we describe hom-connections on a quasi-free (or smooth) algebra, hom-connections with respect to differential graded algebras based on derivations, and integral forms on the matrix algebra and on the quantum plane. 

\begin{example}
Assume that $A$ is a finitely generated algebra. $A$ is said to be {\em quasi-free} \cite{CunQui:alg} or {\em smooth} \cite{Sch:smo}, provided that the universal one-forms $\oaun 1$ are (finitely generated) projective as an $A$-bimodule. Let $\omega_i\in \oaun 1$, ${\zeta}_i \in \lrhom AA{\oaun 1}{A\ot A}$ be a finite dual basis. Define
$$
\bar{\zeta}_i= \mu\circ {\zeta}_i \in  \lrhom AA{\oaun 1}{A}. 
$$
The $A$-bilinearity of the $\bar{\zeta}_i$ implies that the maps
$$
\partial_i: A\longrightarrow  A, \quad  a \longmapsto \bar{\zeta}_i(1\ot a - a\ot 1),
$$
are derivations (not twisted, i.e.\ each $\partial_i$ is a $q$-skew derivation with $q=1$ and $\sigma_i = \id$). Let 
$$
\oan 1: = \oaun 1/\left[A, \oaun 1\right],
$$
(note that $\oan 1 = H_0(A,\oaun 1)$, the zeroth Hochschild homology of $A$ with values in $\oaun 1$). Write $\pi: \oaun 1\to \oan 1$ for the canonical $A$-bimodule epimorphism and then define
$$
d: A\to \oan 1, \qquad a\mapsto \pi(1\ot a - a\ot 1), 
$$
and $ \oomega_i = \pi(\omega_i)$. Since, in the universal calculus, $da = \sum_i {\zeta}_i(da) \omega_i$ (where we view the $A$-bimodule $\oan 1$ as a left module over the enveloping algebra $A^e$), a straightforward calculation yields 
$$
da = \sum_i \partial_i(a)\oomega_i,
$$
in $\oan 1$. By construction, $\oan 1$ is a central $A$-bimodule, hence $\oomega_i a = a\oomega_i$. Thus, by Corollary~\ref{cor.main}, a quasi-free algebra $A$ admits  a hom-connection   $\nabla: \rhom A {\oan 1} A\to A$, 
$$
\nabla(f) = \sum_i \partial_i\left(f\left(\oomega_i\right)\right) = \sum_i \mu\circ \zeta_i\left ( \underset{}{} 1\ot f\left(\pi\left(\omega_i\right)\right) - f\left(\pi\left(\omega_i\right)\right)\ot 1\right) .
$$ 
\end{example}

\begin{example}\label{ex.diff}
The construction described in Theorem~\ref{thm.hom-der} is also applicable to differential graded algebras based on derivations introduced in \cite{Dub:der}; see \cite{Dub:lec} for a review and various applications e.g.\  to Yang-Mills theories. Let $A$ be an algebra and set $R = Z(A)$ (the centre of $A$). Denote by $D(A)$ the Lie algebra of all derivations $A\to A$. Take a Lie subalgebra and $R$-submodule $\fg \subseteq D(A)$, and define $\oan n$ as a set of $R$-multilinear antisymmetric maps $\fg^{\ot_R n} \to A$. $\oa = \oplus_i \oan i$ is an algebra with the product
$$
\omega\eta(\cX_1,\ldots , \cX_{p+q}) = \frac{1}{p!q!} \sum_{\pi \in S_{p+q}} (-1)^{\mathrm{sgn}(\pi)}\omega(\cX_{\pi(1)}, \ldots , \cX_{\pi(p)})\eta (\cX_{\pi(p+1)}, \ldots , \cX_{\pi(p+q)}). 
$$
The differential is given by the Koszul formula
\begin{eqnarray*}
d(\omega)(\cX_1, \ldots , \cX_{n+1}) &=& \sum_{i=1}^{n+1} (-1)^{i+1} \cX_i\left(\omega(\cX_1, \ldots, \cX_{i-1}, \cX_{i+1}, \ldots,  \cX_{n+1})\right)\\
&& \hspace{-1in} + \sum_{1\leq i<j\leq n+1} (-1)^{i+j} \omega([\cX_i,\cX_j], \cX_1, \ldots, \cX_{i-1}, \cX_{i+1}, \ldots, \cX_{j-1}, \cX_{j+1}, \ldots , \cX_{n+1}).
\end{eqnarray*}
In particular $da(\cX) = \cX(a)$. Suppose that $\fg$ is finitely generated and projective as a right $R$-module, and let $\cX_1, \ldots, \cX_n \in \fg$, $\oomega_1, \ldots, \oomega_n \in \fg^* = \rhom R {\fg} R$ be a dual basis. Then $\oan 1 \simeq \fg^*\ot_R A$ as $A$-bimodules. The isomorphism is 
$$
\vartheta :\oan 1 \to \fg^*\ot_R A, \qquad \omega \mapsto \sum_i \oomega_i \ot_R \omega(\cX_i).
$$
$\fg^*\ot_R A$ is a left $A$-module by $a(\oomega \ot_R b) = \oomega \ot_R ab$. The differential $d: A\to \oan 1$ translates through $\vartheta$ to 
\begin{eqnarray*}
da &:=& \vartheta(da) = \sum_i \oomega_i\ot_R da(\cX_i) = \sum_i \oomega_i \ot_R \cX_i(a)\\
&=& \sum_i \cX_i(a)(\oomega_i\ot_R 1) = \sum_i (\oomega_i\ot_R 1)\cX_i(a).
\end{eqnarray*}
Thus $\oan 1 \simeq \fg^*\ot_R A$ is generated (as a left and right $A$-module) by
$
\omega_i := \oomega_i \ot_R 1.
$
Furthermore, for all $a\in A$,
$$
a\omega_i = \omega_i a, \qquad da = \sum_i \cX_i(a)\omega_i.
$$
Hence there is a hom-connection 
\begin{equation}\label{nablaD1}
 \nabla(f) \,=\, \sum_i \cX_i(f(\omega_i)), \quad \text{for all } f \in \rhom A {\oan 1} A.
\end{equation}
In the light of the chain of the ($A$-bimodule) isomorphisms
$$
 \rhom A {\oan 1} A \simeq \rhom R {\fg^*} A \simeq A\ot_R\fg \simeq \fg\ot_R A,
$$
the formula for the hom-connection comes out as
\begin{equation}\label{nablaD2}
\nabla(\cX \otimes_R a) \,=\, \sum_i \cX_i(\oomega_i(\cX)a), \quad \text{for all } \cX \in \fg, a \in A.
\end{equation}
In case $\fg = D(A)$, this last formula is a non-commutative version of the example of a co-connection constructed in \cite[Section 3, Example]{Vin:coc}.

Assume now that $\fg$ is free as an $R$-module with a finite basis $\cX_i$, and let $\oomega_i$ be the dual basis, i.e.\ $\oomega_i(\cX_j) = \delta_{ij}$. Set $\omega_i = \oomega_i \ot_R 1_A$ as before. Since $\fg$ is a Lie-subalgebra of $D(A)$ generated as an $R$-module by the $\cX_i$, there are elements $c_{ijl} \in R$ such that
\begin{equation}\label{exes}
[\cX_i,\cX_j] = \sum_{l}c_{ijl}\cX_l.
\end{equation}
Then one finds that 
\begin{equation}\label{omegas}
\omega_i\, \omega_j = -\omega_j\, \omega_i, \qquad d\omega_l = -\frac{1}{2} \sum_{i,\, j} c_{ijl}\omega_i\, \omega_j.
\end{equation}
Using these expressions and the derivation property of each of the $\cX_i$, the curvature $F$ of hom-connection \eqref{nablaD1} can be computed as, for all $f\in \rhom A {\oan 2} A$, 
$$
F(f) = - \frac{1}{2} \sum_{i,\, j, \, l} f\left(\cX_l\left(c_{ijl}\right)\omega_i\, \omega_j\right).
$$
In particular, if $R=\k $, then $\nabla$ in \eqref{nablaD1} is a flat hom-connection.
\end{example}

\begin{example}\label{ex.mn}
As a special case of Example~\ref{ex.diff}, take $\k  =\C$, $A = M_n(\C)$ (the algebra of complex $n\times n$-matrices) and $\fg = D(A) = sl(n,\C)$ (the Lie algebra of complex traceless $n\times n$-matrices). This is an example of non-commutative geometry studied in \cite{DubKer:non}. In this case $R=Z(A)=\C$, and thus the constructed hom-connection is flat. To get further insight into this example, choose a basis $E_l$ for $sl(n,\C)$ (e.g.\  a basis consisting of Hermitian matrices). Then the corresponding basis  $\{\cX_l\; ,  l=1,\ldots , n^2-1\}$ for $D(A)$ can be chosen as $\cX_l(a) = \imath [E_l, a]$. The formula \eqref{nablaD2} then yields a flat hom-connection
$$
\nabla\left(\sum_{l} \cX_l \ot a_l \right) = \sum_l \cX_l (a_l) = \sum_l \imath  [E_l, a_l].
$$
This last expression affirms that $\im (\nabla) = sl(n,\C)$, therefore $\coker (\nabla) = \C$, and the $\nabla$-integral $\Lambda: M_n(\C) \to \C$ comes out as
$$
\Lambda (a) = \Lambda\left(\left(a-\frac{1}{n} \tr (a)\right)+\frac{1}{n} \tr (a)\right) = \frac{1}{n} \tr (a),  
$$
since, by definition, $\Lambda$ vanishes on all traceless matrices (the image of $\nabla$). This is exactly the integral on $M_n(\C)$ considered in \cite[Section~VA]{DubKer:non}. 

One can construct an isomorphism between the de Rham complex and the complex of integral forms as follows. Set $N=n^2-1$ and suppose that the matrices $E_l$, $l=1,\ldots, N$ form a fundamental representation of $sl(n,\C)$, i.e.\ the corresponding structure constants $c_{ijk}$ in equations \eqref{exes}--\eqref{omegas} are completely antisymmetric. Consider the following diagram 
$$
\xymatrix{ A\ar[r]^d \ar[d]_{\Phi_0} & \oan 1\ar[r]^d \ar[d]_{\Phi_1} & \oan 2\ar[r]^-d\ar[d]_{\Phi_2} & \ldots  \ar[r]^-d & \oan {N-1} \ar[r]^d\ar[d]_{\Phi_{N-1}} &\oan N\ar[d]^{\Phi_N}\\
\oan N^* \ar[r]^{\nabla_{N-1}} & \oan {N-1}^*\ar[r]^{\nabla_{N-2}} & \oan {N-2}^*\ar[r]^-{~\nabla_{N-3}} & \ldots  \ar[r]^-{\nabla_1} &\oan 1^*\ar[r]^\nabla & A  \, ,}
$$
where $\Phi_N$ is the canonical isomorphism given by 
$$
\Phi_N(\omega_1\, \omega_2 \cdots \omega_N) = 1,
$$
and, for all $\omega \in \oan k$, 
$$
\Phi_k(\omega): \oan {N-k} \longrightarrow A, \qquad \omega' \longmapsto  (-1)^{(N-1)k}\Phi_N(\omega \, \omega') .
$$
Since the structure constants $c_{ijk}$ are non-zero only when all the indices are different, equations \eqref{omegas} imply that
$$
d(\omega_1\, \omega_2\, \cdots \omega_{k-1}\, \omega_{k+1} \cdots \omega_N) = 0,
$$
for all $k=1,\ldots ,N$. Using this and the definition of the hom-connection $\nabla$ one easily checks that the right-most square in the above diagram is commutative. Combining the commutativity of the right-most square with the definitions of $\nabla_k$ in equations  \eqref{nablan} one can check that all the remaining squares are commutative. Since $\Phi_N$ is an isomorphism (of $A$-bimodules), each of the $\Phi_k$ is an injective map.  A simple dimension-counting argument then yields that the $\Phi_k$  are isomorphisms.
\end{example}

\begin{example}\label{ex.q.plane}
This example deals with a (quantum-group) covariant differential calculus on the quantum hyperplane introduced in \cite{PusWor:twi}, \cite{WesZum:cov}. Although derived in the context of quantum groups, and hence really belonging to forthcoming Section~\ref{sec.examples}, the knowledge of this context is not needed here. To make the presentation more succinct we discuss only the case of the two-dimensional quantum plane, but with a two-parameter differential structure. 

The quantum plane is a unital algebra $A$ generated by $x,y$ subject to relation $xy=qyx$, where $q$ is a non-zero element of $\k $, i.e.\ $A=\k [x,y]/\langle xy-qyx\rangle$. $\oan 1$ is generated freely by $\omega_1 = dx$ and $\omega_2 =dy$, subject to relations
\begin{equation}\label{q.plane.one-forms}
dx x = p xdx, \quad dx y = qydx + (p-1)xdy, \quad dy x = pq^{-1}xdy, \quad dy y = p ydy , 
\end{equation}
where $p$ is a non-zero elelment of $\k $; see \cite[Theorem~2.1]{PusWor:twi} or \cite[Section~2A]{BrzDab:dif}. This first order differential calculus extends to the full differential graded algebra $\oa = A \oplus \oan 1 \oplus \oan 2$, in which  $dydx = -pq^{-1}dxdy$ and  $(dx)^2 = (dy)^2 =0$. The corresponding matrix $\sigma$ is
$$
\sigma(x^ry^s) = \begin{pmatrix} 
p^rq^{s}x^ry^s & p^{r}(p^{s} -1)x^{r+1}y^{s-1} \\
0 & p^{r+s}q^{-r}x^ry^s
\end{pmatrix}.
$$
Since $\sigma$ is upper-triangular and its diagonal entries are bijective, the corresponding twisted multi-derivation is free. The construction in the proof of Proposition~\ref{triangular} yields
$$
\bsi(x^ry^s) = \begin{pmatrix} 
p^{-r}q^{-s}x^ry^s & 0 \\
p^{-r}q^{r-s+1}(p^{-s} -1)x^{r+1}y^{s-1} & p^{-r-s}q^{r}x^ry^s
\end{pmatrix}, 
$$
and
$$
\hsi(x^ry^s) = \begin{pmatrix} 
p^rq^{s}x^ry^s & p^{r+1}(p^{s} -1)x^{r+1}y^{s-1} \\
0 & p^{r+s}q^{-r}x^ry^s
\end{pmatrix}.
$$
One can now construct a hom-connection on $A$ as in Theorem~\ref{thm.hom-der}. We  concentrate on the following questions: is this hom-connection flat, what is the form of the associated integral and how is  the complex of integral forms related to the de Rham complex. 

Write $\xi_x\in \rhom A {\oan 1} A$ for the dual of $dx$, $\xi_y\in  \rhom A {\oan 1} A$ for the dual of $dy$, and $\xi\in  \rhom A {\oan 2} A$ for the dual of $dxdy$. Then $\xi dx = \xi_y$ and $\xi dy = -pq^{-1}\xi_x$. Since the hom-connection $\nabla$ constructed in Theorem~\ref{thm.hom-der} has the property $\nabla(\xi_x) = \nabla(\xi_y) = 0$, one immediately concludes that $\nabla_1(\xi) = 0$. Since every element of $ \rhom A {\oan 2} A$ is of the form $\xi a$, for some $a\in A$, and the curvature $F$ of $\nabla$ is $A$-linear, we conclude that $F=0$. 

Noting that, for all $a\in A$, $\nabla(\xi_y a) = \xi_y(da)$ one easily finds that
$$
x^ry^s = p^{r+s}q^{-r} \frac{p-1}{p^{s+1}-1}\nabla(\xi_y x^ry^{s+1}). 
$$
Therefore, $\nabla: \rhom A {\oan 1} A \to A$ is an epimorphism, so $\coker (\nabla) = 0$, and thus the $\nabla$-integral is zero. 

Finally, one easily checks the commutativity of  the following diagram
$$
\xymatrix{ A\ar[rr]^d \ar[d]_{\Theta^*} && \oan 1\ar[rr]^d \ar[d]_\Phi &&\oan 2\\
\oan 2^* \ar[rr]^{\nabla_1} && \oan 1^*\ar[rr]^{\nabla} &&  A \ar[u]_\Theta \, ,}
$$
where, for all $a,b\in A$,  
$$
\Theta(a) = dxdy a, \qquad \Theta^*(a) = \xi a, \qquad \Phi(dx a + dy b) = pq^{-1}\xi_x b - \xi_ya . 
$$
All the vertical arrows are (right $A$-module) isomorphisms. Consequently, the integral complex is isomorphic to the de Rham complex. 
\end{example}

\section{Examples of integral forms on quantum groups and spaces}\label{sec.examples}
\setcounter{equation}{0}

The data which enter the definition of a differential calculus, and hence also the  hom-connection  in Theorem~\ref{thm.hom-der}, have natural interpretation in terms of actions of coalgebras and Hopf algebras. 
First, there is a bijective correspondence between algebra maps $\sigma:A\rightarrow M_n(A)$ and {\em measurings} of the $n\times n$ comatrix coalgebra $M_n^c(\k )$ to $A$, and hence right (or left) module algebra structures of $A$ over $\mathcal{O}(M_n(\k ))$. Write $\theta_{ij}$ for a  basis for $M_n^c(\k )$ with comultiplication $\Delta(\theta_{ij})=\sum_k \theta_{ik}\otimes\theta_{kj}$ and counit $\eps(\theta_{ij})=\delta_{ij}$.  If $H$ is a Hopf algebra containing $M_n^c(\k )$, e.g.\ $H=\mathcal{O}(SL(n))$, such that $A$ is a right $H$-module algebra, then the associated algebra map $\sigma$ is  $\sigma_{ij}(a):=a\ra \theta_{ij}$. Furthermore, assignments $\bsi_{ij}(a)=a\ra S(\theta_{ji})$ and $\hsi_{ij}(a)=a\ra S^2(\theta_{ij})$, where $S$ is the antipode of $H$, define maps $\bsi$ and $\hsi$ as in Definition~\ref{def.der.free}. If $A$ is a left $H$-module algebra and $H$ has a bijective antipode, then $\sigma$, $\bsi$, $\hsi$ are given by $\sigma_{ij}(a) = \theta_{ij}\la a$, $
 \bsi_{ij}(a)= S^{-1}(\theta_{ji})\la a$, $\hsi_{ij}(a)= S^{-2}(\theta_{ji})\la a$; see Theorem~\ref{thm.cov} below for the proof. 

Extending the comatrix coalgebra by an $n+1$-dimensional vector space $V$ with basis $g,\chi_1,\ldots, \chi_n$ we define a coalgebra $C_n=M_n^c(\k )\oplus V$ with comultiplication $\Delta{(g)}=g\otimes g$ and $\Delta(\chi_i)= g\otimes \chi_i + \sum_j \chi_j\otimes \theta_{ji}$ and counit $\eps(g)=1$ and $\eps(\chi_i)=0$. Then right twisted $n$-multi-derivations $(\partial, \sigma)$ on $A$ correspond bijectively to  measurings of $C_n$ to $A$. 

These observations are a basis for finding examples of twisted multi-derivations, and so of hom-connections. We describe  these examples presently from a differential geometric point of view.

\subsection{Quantum groups with (left) covariant differential calculi}\label{sec.sl2}
Let $A$ be a Hopf algebra with the coproduct $\Delta$, counit $\eps$ and bijective antipode $S$. Following \cite{Wor:dif}, a first order differential calculus $\oan 1$ on a
Hopf algebra $A$ is said to be {\em left covariant}, if the
coproduct $\Delta$ extends to a map $\Delta_L :\oan 1 \to
A\otimes \oan 1$ by the formula
$$
\Delta_L(ad(b)) = a\sw 1b\sw 1\otimes a\sw 2 d (b\sw 2),
$$
where the Sweedler notation $\Delta(a) = a\sw 1\ot a\sw 2$ (summation implicit) is used. By \cite[Theorem~2.1, Theorem~5.2]{Wor:dif}  the whole information about a left covariant differential calculus is contained in the following data: elements $\omega_i\in \oan 1$ such that $\Delta_L(\omega_i) = 1\ot \omega_i$ (one says $\omega_i$ are {\em left-invariant}) and that freely generate $\oan 1$ as an $A$-module, and $\k $-linear maps $\theta_{ij}, \chi_i: A\to \k $, $i,j=1,2,\ldots ,n$. These satisfy the following relations, for all $a,b\in A$,
\begin{equation} \label{eq.theta}
\theta_{ij}(ab) = \sum _k \theta_{ik}(a)\theta_{kj}(b), \qquad \theta_{ij}(1) = \delta_{ij},
\end{equation}
\begin{equation}\label{eq.chi}
\chi_i(ab) = \sum _j \chi_j(a) \theta_{ji}(b) + \eps(a)\chi_i(b).
\end{equation}
$A^* = \rhom \k  A \k $ is an algebra with convolution product  $f*g = (f\ot g)\circ \Delta$. $A$ is a left $A^*$-module with the multiplication $f\la a = (\id\ot f)(\Delta(a))$.  Using this notation the commutation rules in $\oan 1$ and the definition of the exterior differential $d: A\to \oan 1$ are given by
\begin{equation}\label{eq.omega.cov}
\omega_i a = \sum_ j(\theta_{ij}\la a)\omega_j, \qquad d(a) = \sum _i (\chi_i\la a)\omega_i.
\end{equation}
Note that relations \eqref{eq.theta}--\eqref{eq.chi} mean that the elements $\theta_{ij}$, $\chi_i$ span a subcoalgebra in the Hopf dual $A^\circ\subseteq A^*$  of $A$ with coproduct and counit (in $A^\circ$), $\Delta(\theta_{ij})=\sum_k \theta_{ik}\otimes\theta_{kj}$,  $\Delta(\chi_i)= 1\otimes \chi_i + \sum_j \chi_j\otimes \theta_{ji}$,  $\eps(\theta_{ij})=\delta_{ij}$  and $\eps(\chi_i)=0$ (recall that the unit in $A^*$ is given by the counit in $A$). Since $A$ is a left $A^\circ$-module algebra (with multiplication $\la$) the observations made in the preamble to Section~\ref{sec.examples} yield the following
\begin{theorem}\label{thm.cov}
Let $\oan 1$ be a left covariant differential calculus on a Hopf algebra $A$ with bijective antipode $S$. 
\begin{zlist}
\item  Let $\{\omega_1,\omega_2, \ldots, \omega_n\}$ be a left invariant basis for $\oan 1$,  let $\{\xi_1,\xi_2, \ldots, \xi_n\}$ be  its  right dual basis and let $\{\lambda_1,\lambda_2, \ldots, \lambda_n\}$ be its left dual basis. Then
$$
\xymatrix@R=0pt{ \nabla : \rhom A {\oan 1} A \ar@{->}[r]  & A \\ f \ar@{|->}[r]  &  \sum_i \left(S^2 \circ \lambda_i\circ d \circ S^{-2}\right) \left(  f(\omega_i)\right) 
= \sum_i \left(\chi_i\circ S^{-2}\right) \la  f(\omega_i),  } 
$$
where $\chi_i = \eps\circ \lambda_i\circ d: A\to \k $, is  a unique hom-connection  on $A$ such that $\nabla (\xi_i)=0$, for all $i=1,2,\ldots ,n$.
\item Assume that $\oan 1$ extends to a differential graded algebra such that the hom-connection $\nabla$ constructed in item (1) is flat, and that there is  a right integral $\lambda: A\to \k $ on the Hopf algebra $A$. Then there exists a unique map $\varphi: \coker(\nabla) \to \k $ such that
$$
\lambda = \varphi \circ\Lambda ,
$$
where $\Lambda: A \to  \coker(\nabla)$ is the $\nabla$-integral on $A$.
\end{zlist}
\end{theorem}
\begin{proof}
 (1) Note that the $\chi_i$  are the same as in \eqref{eq.omega.cov} and let $\theta_{ij}$ be the corresponding (as in \eqref{eq.omega.cov}) maps $A\to \k $. The following defines a free right twisted multi-derivation $(\partial,\sigma;\bsi, \hsi)$,
$$
\partial_i (a) = \chi_i \la a, \quad \sigma_{ij} (a) = \theta_{ij}\la a, \quad \bsi_{ij}(a) = (\theta_{ji}\circ S^{-1})\la a, \quad  \hsi_{ij}(a) = (\theta_{ij}\circ S^{-2})\la a.
$$
 Indeed, first note that, since the antipode is an anti-algebra map, comultiplication is an algebra map and $\theta_{ij}$ satisfy conditions \eqref{eq.theta}, all three  maps $\sigma$, $\hsi$, $\bsi$ are algebra morphisms (note the distribution of indices in the definition of $\bsi$). Equations \eqref{eq.chi} force $(\partial, \sigma)$ to be a right twisted multi-derivation. Checking that relations   \eqref{bar.sigma} and \eqref{hat.sigma} hold is performed by  the standard gymnastics with the Sweedler notation. For example, to prove the first of \eqref{hat.sigma}, take any $a\in A$ and compute
\begin{eqnarray*}
\sum_k \hsi_{ik}\circ \bsi_{jk}(a) &=& \sum_k a\sw 1\theta_{ik} (S^{-2}( a\sw 2))\theta_{kj}(S^{-1} (a\sw 3)) = a\sw 1 \theta_{ij}(S^{-2}(a\sw 2)S^{-1}(a\sw 3)) \\
 &=& a\sw 1 \theta_{ij}(S^{-2}(a\sw 2 S(a\sw 3))) = a\theta_{ij}(1) = a\delta_{ij},
 \end{eqnarray*}
 where the second and last equalities follow by  \eqref{eq.theta}, and the third and fourth one use properties of the antipode. The remaining equations \eqref{bar.sigma} and \eqref{hat.sigma}  are checked in the same manner. 
 
Even the most perfunctory comparison of equations \eqref{eq.omega.cov} with \eqref{eq.dif.rel} and \eqref{eq.dif.d} reveals  that the left covariant differential calculus we start with (and which is determined by the $\chi_i$, $\theta_{ij}$)  is the same as the differential calculus constructed from the right twisted multi-derivation $(\partial,\sigma;\bsi, \hsi)$.  Theorem~\ref{thm.hom-der} implies the existence and uniqueness of a hom-connection $(A,\nabla )$ such that $\nabla (\xi_i) =0$. The formula for $\nabla $ given in the proof of Theorem~\ref{thm.hom-der} can be simplified as follows:
\begin{eqnarray*}
\nabla ( f) \!\! &=& \!\! \sum_i \partial_i^\sigma ( f(\omega_i)) = \sum_{i,j,k} f(\omega_i)\sw 1\theta_{jk}(S^{-1}(f(\omega_i)\sw 2)) \chi_j( f(\omega_i)\sw 3)\theta_{ki}(S^{-2}(f(\omega_i)\sw 4))\\
&\!\!\!\! \!\!\!\!\!\!=&\!\!\!\!\!\!\!\!   \sum_{i,j}  f(\omega_i)\sw 1 \chi_j( f(\omega_i)\sw 3)\theta_{ji}(S^{-1}(f(\omega_i)\sw 2)S^{-2}(f(\omega_i)\sw 4))\\
&\!\!\!\! \!\!\!\!\!\!=&\!\!\!\!\!\!\!\!  \sum _i f(\omega_i)\sw 1 \left[ \underset{}{} \chi_i( f(\omega_i)\sw 3 S^{-1}(f(\omega_i)\sw 2)S^{-2}(f(\omega_i)\sw 4)) \right.\\
&& \hspace{1in} -  \left. \underset{}{} \eps(f(\omega_i)\sw 3)\chi_i(S^{-1}(f(\omega_i)\sw 2)S^{-2}(f(\omega_i)\sw 4)) \right]\\
&\!\!\!\!\!\!\!\! \!\!=&\!\!\!\!\!\!\!\!  \sum_i f(\omega_i)\sw 1\chi_i(S^{-2}(f(\omega_i)\sw 2)) = \sum_i (\chi_i\circ S^{-2}) \la  f(\omega_i),
\end{eqnarray*}
where the third equality follows by \eqref{eq.theta} and the fourth  by \eqref{eq.chi}. The fifth equality follows by the properties of the antipode and counit and by  \eqref{eq.chi} (to conclude that $\chi_i(1) =0$). Thus the unique hom-connection on $A$ such that $\nabla (\xi_i) =0$ has the form stated. 

(2) Recall that a right integral on a Hopf algebra $A$ is a $\k $-linear map $\lambda: A\to \k $ such that, for all $a\in A$,
\begin{equation}\label{eq.integral}
\lambda(a\sw 1)a\sw 2 = \lambda(a),
\end{equation}
i.e.\ $\lambda$ is a right $A$-colinear map. For all $f\in  \oan 1^*$,
$$
\lambda\left(\nabla\left(f\right)\right) 
= \sum_i \lambda(f(\omega_i)\sw 1)\chi_i(S^{-2}(f(\omega_i)\sw 2) ) = 
\sum_i \chi_i \circ S^{-2}\left(\lambda\left(f(\omega_i)\right) 1 \right) =0,
$$
where the second  equality follows by \eqref{eq.integral} and the final equality is a consequence of $\chi_i(1)=0$.
Therefore, $\lambda\circ\nabla = 0$. By the universality of cokernels there exists a unique $\k $-linear map $\varphi: \coker(\nabla) \to  \k $ completing the following diagram
$$
\xymatrix{ \oan 1^* \ar[rr]^{ \nabla} && A\ar[rr]^-{\Lambda } \ar[d]_{\lambda} &&  \coker(\nabla)\ar@{-->}[dll]^\varphi\\
&& \k  . && }
$$
This completes the proof.
\end{proof}

The 3-dimensional (or 3D) calculus on the quantum group $\cO_q(SL(2))$ defined in \cite{Wor:twi} provides an example of the construction described in Theorem~\ref{thm.cov}. Suppose that $\k $ is a field of characteristic 0 (typically $\k =\C$). $\cO_q(SL(2))$ is a Hopf algebra generated by
$\{\alpha,\beta,\gamma,\delta\}$ with the relations
\begin{eqnarray}\label{comm.sl2}
&& \alpha\beta\,=\, q\,\beta\alpha\ ,\quad
\alpha\gamma\,=\, q\,\gamma\alpha\ ,\quad
\beta\gamma\,=\, \gamma\beta\ ,\quad
\beta\delta\,=\, q\,\delta\beta\ ,\quad
\gamma\delta\,=\, q\,\delta\gamma\ , \cr
&&\alpha\delta\,=\, \delta\alpha+(q-q^{-1})\,\beta\gamma\ ,\quad
\alpha\delta-q\,\beta\gamma\,=\,1\ ,
\end{eqnarray}
where $q\neq 0$ is a scalar which is not a root of unity. 
The coproduct is given by
\begin{eqnarray}\label{cop.sl2}
&&\Delta(\alpha)\,=\,\alpha\tens\alpha+\beta\tens\gamma\, ,\quad
\Delta(\beta)\,=\,\alpha\tens\beta+\beta\tens\delta\, ,\cr
&& \Delta(\gamma)\,=\,\gamma\tens\alpha+\delta\tens\gamma\, ,\quad
\Delta(\delta)\,=\, \delta\tens\delta+\gamma\tens\beta\, ,
\end{eqnarray}
and counit and antipode are
\begin{eqnarray*}
&& \eps(\alpha)\,=\,\eps(\delta)\,=\,1\ ,\quad
\eps(\beta)\,=\,\eps(\gamma)\,=\,0\ ,\cr
&& S(\alpha)\,=\,\delta\ ,\quad
S(\delta)\,=\,\alpha\ ,\quad
S(\beta)\,=\,-q^{-1}\,\beta\ ,\quad
S(\gamma)\,=\,-q\,\gamma\ .
\end{eqnarray*}
$\cO_q(SL(2))$ is a $\Z$-graded algebra with the grading defined on generators by $|\alpha| = |\gamma| = 1$, $|\beta| = |\delta| = -1$. If $\k = \C$ and $q\in \R$, $\cO_q(SL(2))$ can be equipped with a $*$-Hopf algebra structure with  $\alpha^*=\delta$ and $\beta^* = -q\gamma$, thus giving rise to the real form $\cO_q(SU(2))$ of $\cO_q(SL(2))$. If $q\in (0,1)$, then $\cO_q(SU(2))$ is a dense subalgebra of the $C^*$-algebra of continuous functions on the quantum group $SU_q(2)$. The forthcoming purely algebraic discussion applies to this topological situation too.

The 3D left covariant differential calculus on $A= \cO_q(SL(2))$ 
is generated by three left invariant one-forms
$\{\omega_0,\omega_+,\omega_-\}$ that are required to satisfy the following commutation relations
\begin{eqnarray}\label{eq.3Drel}
\omega_0\,\alpha \,=\, q^{-2}\,\alpha\,\omega_0\, , &&
\omega_0\,\beta \,=\, q^2\,\beta\,\omega_0\ ,\cr
\omega_+\,\alpha \,=\, q^{-1}\,\alpha\,\omega_+\, , &&
\omega_+\,\beta \,=\, q\,\beta\,\omega_+\ ,\cr
\omega_-\,\alpha \,=\, q^{-1}\,\alpha\,\omega_-\, , &&
\omega_-\,\beta \,=\, q\,\beta\,\omega_-\, ,
\end{eqnarray}
and similarly for replacing $\alpha\to\gamma$ and $\beta\to\delta$. The action of the exterior differential $d$ on the generators 
is 
\begin{eqnarray} \label{eq.3Drel.d}
d (\alpha) \,=\, \alpha\,\omega_0-q\,\beta\,\omega_+\, ,\quad
d (\beta) \,=\, -q^2\,\beta\,\omega_0 + \alpha\,\omega_-\, ,\cr
d (\gamma) \,=\, \gamma\,\omega_0-q\,\delta\,\omega_+\, ,\quad
d (\delta) \,=\, -q^2\,\delta\,\omega_0+ \gamma\,\omega_-\, .
\end{eqnarray}

The form of relations \eqref{eq.3Drel} immediately reveals that the matrix $\sigma$ is diagonal with the diagonal entries $\sigma_0$ and $\sigma_+=\sigma_-$ given by, for all homogeneous $a\in A$ (with the $\Z$-degree $|a|$),
\begin{equation}\label{sigma.3d}
\sigma_0(a) = q^{-2|a|}a, \qquad \sigma_\pm(a) = q^{-|a|}a.
\end{equation}
Equations \eqref{eq.3Drel.d} determine $\sigma_i$-twisted derivations $\partial_i$. Explicitly, in terms of actions on generators of $A$ these are
\begin{eqnarray} \label{eq.partial.pm}
&&\partial_0(\alpha) = \alpha,  \quad   \partial_0(\beta) = -q^2\beta,  \quad  \partial_0(\gamma) = \gamma,   \quad  \partial_0(\delta) = -q^2\delta \, ,\nonumber \\
&& \partial_+(\alpha) = -q\beta, \quad  \partial_+(\beta) = 0,  \quad  \partial_+(\gamma) = -q\delta,  \quad \partial_+(\delta) = 0 \, ,\\
&& \partial_-(\alpha) = 0, \quad  \partial_-(\beta) = \alpha, \quad  \partial_-(\gamma) = 0, \quad  \partial_-(\delta) = \gamma \, . \nonumber
\end{eqnarray}
The maps  $\partial_0, \partial_+,\partial_-$ are $q_i$-skew derivation with constants 1, $q^{-2}$ and $q^2$. respectively. Therefore, Theorem~\ref{thm.cov} or Theorem~\ref{thm.hom-der} give rise to the following hom-connection on $\cO_q(SL(2))$,
\begin{equation}\label{eq.hom.3d}
\nabla  (f) = \partial_0\left( f\left(\omega_0\right)\right) + q^{-2}\partial_+\left( f\left(\omega_+\right)\right) + q^{2}\partial_-\left( f\left(\omega_-\right)\right).
\end{equation}

In \cite{Wor:twi} Woronowicz describes the full differential graded algebra built on the 3D calculus. The relations for  the higher forms are 
\begin{equation}\label{eq.3Drel.higher}
\omega_i^2  =  0, \quad 
 \omega_+\omega_-=-q^2\,\omega_-\omega_+\, ,\quad
\omega_0\omega_-=-q^4\,\omega_-\omega_0\, ,\quad
\omega_+\omega_0\,=\,-q^4\,\omega_0\omega_+\, ,
\end{equation}
and the exterior derivative is
\begin{equation}\label{eq.3Drel.d.higher}
d (\omega_0) =q\omega_-\omega_+\, ,\quad
d (\omega_+) =q^2(q^2+1)\,\omega_0\omega_+\, , \quad 
d (\omega_-) =q^2(q^2+1)\omega_-\omega_0\, .\quad
\end{equation}
In degree $3$, $\Omega^3(A)$ is generated by the (volume) form $\omega_-\omega_0\omega_+$. 

\begin{proposition}\label{prop.3d.flat}
Let $A=\cO_q(SL(2))$ and $\oa$ be the differential graded algebra corresponding to the 3D calculus (and described above). The hom-connection \eqref{eq.hom.3d} is flat.
The associated complex of integral forms 
$$
\xymatrix{
\oan 3^* \ar[rr]^{\nabla_2} && \oan 2^*\ar[rr]^{\nabla_1} && \oan 1^*\ar[rr]^\nabla && A} ,
$$
is isomorphic to the de Rham complex $(\oa , d)$.
\end{proposition}
\begin{proof}
In the light of relations \eqref{eq.3Drel.higher}, the bimodule $\oan 2$ is generated by three forms  $\{ \omega_-\omega_+, \omega_-\omega_0 ,   \omega_0\omega_+\}$. Let $\{\phi_0,\phi_+,\phi_-\}$ be a dual basis, i.e.\ $\phi_i$ are determined by  
$$
\phi_0(\omega_-\omega_+) =1\, ,  \quad \phi_+(\omega_-\omega_0) =1\, , \quad   \phi_-(\omega_0\omega_+) =1\, , 
$$
and zero on other generators. 
Again by inspection of relations \eqref{eq.3Drel.higher} one concludes that any $f\in \rhom A {\oan 2} A$ can be written as
$
f = \phi_{0}a_0 + \phi_+a_+ +\phi_- a_-,$ for suitably defined $a_i\in A$; see the proof of Theorem~\ref{thm.hom-der}. Since the curvature of a hom-connection is a right $A$-linear map, it suffices to compute it on the $\phi_{i}$. One easily computes that
$$
\nabla_1(\phi_0) = q\xi_0, \quad \nabla_1(\phi_+) = q^2(q^2+1)\xi_-, \quad \nabla_1(\phi_-) = q^2(q^2+1)\xi_+, 
$$
where the $\xi_i\in \rhom A {\oan 1} A$ are as in Theorem~\ref{thm.cov}. $\nabla $ given by \eqref{eq.hom.3d} is the unique hom-connection such that $\nabla (\xi_i) =0$, hence $F(\phi_i) =0$, for all $i=1,2,3$, and the hom-connection \eqref{eq.hom.3d} is flat.

Let $\phi\in  \rhom A {\oan 3} A$ denote the map dual to the volume form $\omega_-\omega_0\omega_+$, i.e.\ $\phi(\omega_-\omega_0\omega_+a) = a$. Using the Leibniz rule for hom-connections and the fact that $\nabla(\xi_\pm) = \nabla(\xi_0) = 0$ one easily finds that $\nabla_2(\phi) =0$, and consequently $\nabla_2(\phi a) = \phi d(a)$, for all $a\in A$. Consider the diagram
$$
\xymatrix{ A\ar[rr]^d \ar[d]_{\Theta^*} && \oan 1\ar[rr]^d \ar[d]_\Phi && \oan 2\ar[rr]^d \ar[d]_\Psi &&\oan 3\\
\oan 3^* \ar[rr]^{\nabla_2} && \oan 2^*\ar[rr]^{\nabla_1} && \oan 1^*\ar[rr]^\nabla && A \ar[u]_\Theta \, ,}
$$
in which all columns are (right $A$-module) isomorphisms defined as follows. $\Theta^*(a) = \phi a$, $\Theta (a) = \omega_-\omega_0\omega_+a$, and
$$
\Phi(\omega_-a + \omega_0 b +\omega_+ c) = \phi_- a -q^4\phi_0b + q^6\phi_+c \, ,
$$
$$
\Psi(\omega_-\omega_0a + \omega_-\omega_+ b +\omega_0\omega_+ c) = \xi_+ a -q^4\xi_0b + q^6\xi_-c \, ,
$$
for all $a,b,c\in \cO_q(SL(2))$. The commutativity of this diagram can be checked by a straightforward albeit lengthy calculation.  Therefore, the de Rham and integral complexes are isomorphic as required.
\end{proof}

It is shown in \cite[Section~3]{Wor:twi} that the third de Rham cohomology, $H^3(\cO_q(SL(2)))$, is a one-dimensional $\k $-space, i.e.\ $H^3(\cO_q(SL(2))) =\k $. Furthermore,
the canonical epimorphism $\Omega^3(\cO_q(SL(2))) \to \coker\left(d: \Omega^2(\cO_q(SL(2)))\to \Omega^3(\cO_q(SL(2)))\right) =\k $ can be obtained by composing the normalised integral or the Haar measure on $\cO_q(SL(2))$ with the inverse of the map $\Theta$ defined in the proof of Proposition~\ref{prop.3d.flat}. In the light of Proposition~\ref{prop.3d.flat}, the $\nabla$-integral $\Lambda: A\to H_0(A; A, \nabla)= \k $ on $\cO_q(SL(2))$ is (a scalar multiple of) the normalised integral on the Hopf algebra $\cO_q(SL(2))$. That is
\begin{equation}\label{integral}
\Lambda\left(\left(\beta\gamma\right)^l\right) = (-1)^l \, \frac{q-q^{-1}}{q^{l+1}- q^{-l-1}}, 
\end{equation}
and zero on all other elements $\alpha^k\beta^m\gamma^n$, $\delta^k\beta^m\gamma^n$, $k,m,n\in \N$, $m\neq n$ of the standard linear basis for $\cO_q(SL(2))$; see \cite[Appendix~1]{Wor:com}. 
 
\subsection{Hopf-Galois extensions and the quantum two-sphere}
Let $H$ be a Hopf algebra, and let $A$ be a right $H$-comodule algebra. This means that $A$ is a right $H$-comodule with coaction $\roA: A\to A\ot H$ that is an algebra map (where $A\ot H$ is equipped with the tensor product algebra structure). One often refers to such an $A$ as a {\em quantum space}. 

Similarly to the Hopf algebra case, there is a left action of the convolution algebra $H^*$ on a right $H$-comodule algebra $A$ defined by
$$
f\la a = (\id\ot f)\circ \roA (a), \qquad \forall f\in H^*,\, a\in A.
$$
 A (left) covariant differential calculus $\Omega^1(H)$ on $H$ induces a free right twisted multi-derivation on $A$ as follows. Suppose that  the antipode of $H$ is  bijective, and $\theta_{ij}, \chi_i: H\to \k $ are the data determining $\Omega^1(H)$. Then $(\partial,\sigma ;\bsi, \hsi)$, defined for all $a\in A$ by
 $$
\partial_i (a) = \chi_i \la a, \quad \sigma_{ij} (a) = \theta_{ij}\la a, \quad \bsi_{ij}(a) = (\theta_{ji}\circ S^{-1})\la a, \quad  \hsi_{ij}(a) = (\theta_{ij}\circ S^{-2})\la a,
$$
is a free right twisted multi-derivation on $A$. Thus one can associate a differential graded algebra $\oa$ on $A$ based on $(\partial,\sigma ;\bsi, \hsi)$, and there is a hom-connection $(A,\nabla  :\rhom A {\oan 1} A\to A)$ as in Theorem~\ref{thm.hom-der} or Colloray~\ref{cor.main}. This hom-connection has the same form as the one constructed  in Theorem~\ref{thm.cov}, i.e.\  $\nabla (f)= \sum_i \left(\chi_i\circ S^{-2}\right) \la  f(\omega_i)$, where the $\omega_i$ are generators of $\oan 1$. 

Suppose that a hom-connection $(A,\nabla)$ on a right $H$-comodule algebra $A$ has been constructed. In this section we study the question, when $(A,\nabla)$ descends to a hom-connection on the fixed-point (coinvariant) subalgebra of $B$. In particular, we study the descent of  hom-connections from the total space of a quantum principal bundle to the quantum base space. In algebraic terms, quantum principal bundles are given by  {\em principal comodule algebras}; see e.g.\ \cite{HajKra:pie}. These are examples of {\em Hopf-Galois extensions} whose definition and rudimentary properties we recall presently.

Let $H$ be a Hopf algebra and $A$ a right $H$-comodule algebra, and let  $B=A^{coH} = \{ b\in A \; |\; \roA (b) = b\ot 1\}$ be the subalgebra of coaction invariants (coinvariants). $B\subseteq A$ is called a {\em Hopf-Galois extension}, provided the map
$$
\can: A\ot_B A\to A\ot H, \qquad a\ot_B a'\mapsto a\roA(a'),
$$ 
is bijective. With a Hopf-Galois extension one associates two functors. The {\em induction functor} $-\ot_B A$  sends every right $B$-module $N$ to an $(A,H)$-Hopf module $N\ot_B A$. Recall that an {\em $(A,H)$-Hopf module} or a {\em relative Hopf module} is a vector space $M$ that is a right $A$-module and a right $H$-comodule with coaction $\roM: M\to M\ot H$ such that, for all $a\in A, m\in M$,
$$
\roM(ma) = m\sw 0a\sw 0\ot m\sw 1a\sw 1,
$$
where the Sweedler notation $\roM(m) = m\sw 0\ot m\sw 1$, $\roA(a) = a\sw 0\ot a\sw 1$ for coactions is used. 
The category of $(A,H)$-Hopf modules is denoted by $\M_A^H$. For a right $B$-module $N$, $N\ot_B A$ is a right $A$-module and $H$-comodule by
$$
(n\ot_B a)\cdot a' = n\ot_B aa', \qquad n\ot_B a\mapsto n\ot_B \rho^A(a).
$$
The second functor is the {\em coinvariant functor} $(-)^{coH}$ which sends any object $M\in \M_A^H$ to the right $B$-module
$$
M^{coH} := \{ m\in M \; |\; \roM (m) = m\ot 1\}.
$$
By a theorem of Schneider \cite[Theorem~3.7]{Sch:pri}, the functors $-\ot_B A$ and $(-)^{co H}$ establish an equivalence of categories $\M_A^H$ and $\M_B$, provided $A$ is faithfully flat as a left $B$-module. 

For a Hopf-Galois extension $B\subseteq A$ we choose a {\em covariant} first order differential calculus $\oan 1$. The covariance means that $\oan 1$ is an object in $\M_A^H$ and that $d: A \to \oan 1$ is a right $H$-comodule map. $(\oan 1, d)$ contains the first order calculus on $B$,  $(\obn 1, d)$ (the differential on $B$ is defined by restriction of the differential on $A$). Let
\begin{equation}\label{eq.j}
j: \obn 1\ot_B A\to \obn 1 A\to \oan 1, \qquad \omega\ot_B a\mapsto \omega a. 
\end{equation}
The map $j$ is a morphism in ${}_B\M_A^H$, the category of those $(A,H)$-Hopf modules which are also left $B$-modules by a right $A$-linear right $H$-colinear $B$-action ($(A,\roA)$ is an example).

\begin{theorem}\label{thm.induction}
Let $H$ be a Hopf algebra, and let $B\subseteq A$ be a  Hopf-Galois extension such that $A$ is a faithfully flat left $B$-module. Choose a  covariant  first order differential calculus  $\oan 1$ on $A$ for which there exists a right $H$-colinear right $A$-linear and left $B$-linear map $\Pi: \oan 1 \to \obn 1\ot_B A$ such that $\Pi\circ j = \id$, where $j$ is given by \eqref{eq.j} (i.e.\ $j$ is a section in ${}_B\M_A^H$).  Let $M$ be an $A$-relative $H$-Hopf module. Any hom-connection $(M,\nabla )$ with respect to $\oan 1$ such that 
$$
\nabla  \left( \Rrhom H A {\oan 1} M\right) \subseteq M^{coH},
$$
induces a hom-connection $(M^{coH},\nabla ^{coH})$ with respect to  $\obn 1$.
\end{theorem}
\begin{proof}
Since $B\subseteq A$ is a faithfully flat Hopf-Galois extension, the coinvariant and induction functors are inverse equivalences, and hence there is an isomorphism
$$
\rhom B {\obn 1} {M^{coH}} \ni f \mapsto \widehat{f} \in \Rrhom  HA {\obn 1\ot_B A} M.
$$
Define $\nabla^{coH}: \rhom B {\obn 1} {M^{coH}} \to M^{coH}$, by
 $$ 
 \nabla^{coH}(f) := \nabla(\widehat{f}\circ \Pi ),
 $$
for all $f\in {\rm Hom}_B(\obn 1, M^{coH})$. Note that $\nabla^{coH}(f) \in M^{coH}$ by the assumption on $\nabla$ and the fact that $\widehat{f}\circ \Pi$ is a map in $\M_A^H$. 
Note also that, for all $b\in B$, $\widehat{f}(\Pi(d(b)))=f(d(b))$ since $\Pi{|_{\obn 1}}=\id$. Thus  the $B$-linearity of $\Pi$  and the defining property of $\nabla$ yield 
\begin{eqnarray*}
\nabla^{coH}(f\cdot b) &=& \nabla(\widehat{f\cdot b}\circ \Pi)
= \nabla((\widehat{f}\circ \Pi)\cdot b)\\
&=& \nabla^{coH}(f)b + \widehat{f}(\Pi(d(b))) = \nabla^{coH}(f)b + f(d(b)).
 \end{eqnarray*}
Therefore,  $(M^{coH}, \nabla^{coH})$ is a hom-connection as claimed. 
\end{proof}

In the set-up of Theorem~\ref{thm.induction} suppose that both $(A,\nabla)$ and $(B,\nabla^{coH})$ are flat hom-connections. Let $\Lambda : A\to \coker (\nabla)$ be the $\nabla$-integral and let $\Lambda^{coH}: B\to \coker(\nabla^{coH})$ be the $\nabla^{coH}$-integral. By the construction of $\nabla^{coH}$, $\Lambda |_B \circ \nabla^{coH}  =0$. The universality of cokernels then implies that there exists a unique $\k $-linear map $\varphi: \coker(\nabla^{coH}) \to  \coker(\nabla)$ completing the following diagram
$$
\xymatrix{ \obn 1^* \ar[rr]^{ \nabla^{coH}} && B\ar[rr]^{\Lambda^{coH}} \ar[d]_{\Lambda |_B} &&  \coker(\nabla^{coH})\ar@{-->}[dll]^\varphi\\
&& \coker(\nabla) . && }
$$
This establishes a correspondence between $\nabla$- and $\nabla^{coH}$-integrals.

Although the calculi satisfying requirements of Theorem~\ref{thm.induction} might seem rare, the following corollary asserts the existence of a suitable retraction $\Pi$ for a universal differential calculus on a {\em principal comodule algebra}. In terminology of \cite{HajKra:pie}, a Hopf-Galois extension $B\subseteq A$ by $H$ is said to be a {\em principal comodule algebra} if the antipode of $H$ is bijective and $A$ is a right $H$-equivariantly projective left $B$-module. The latter means that there exists a left $B$-module right $H$-comodule splitting $s: A\to B\ot A$ of the multiplication map $B\ot A\to A$. Such an $s$ can always be normalised so that $s(1) = 1\ot 1$; see \cite{BrzHaj:Che}\footnote{The proofs of claims made in \cite{BrzHaj:Che} are contained in T.\ Brzezi\'nski, P.M.\ Hajac, R.\ Matthes, W.\ Szyma\'nski, {\em The Chern character for principal extensions of noncommutative algebras}, work in progress available at http://www.fuw.edu.pl/ $\widetilde{~~}~$pmh
 /}. As explained in \cite{SchSch:gen}, if the antipode of $H$ is bijective, then $A$ is a principal comodule algebra if and only if it is a faithfully flat (as a left or right $B$-module) Hopf-Galois extension. 

Note that   the universal first order  differential calculus $(\oaun 1 , d)$ on a right $H$-comodule algebra $A$ is $H$-covariant with coaction given by, for all $a,a'\in A$,  
$$
\varrho^{\oaun 1}(d(a)a') = d(a\sw 0)a'\sw 0 \otimes a\sw 1a'\sw 1,
$$ 
i.e.\ $\varrho^{\oaun 1}$ is the restriction of the diagonal coaction of $H$ on $A\ot A$. 

\begin{corollary}\label{lem.princ} Let $H$ be a Hopf algebra with bijective antipode and let $A$ be a principal comodule algebra, $B=A^{coH}$. 
Then $(\obun 1)A$ is a direct summand of $\oaun 1$ in $_B\M_A^H$. 

Consequently, if $(M,\nabla)$ ($M\in \M_A^H)$ is a hom-connection with respect to the universal differential graded algebra, $\oau$, such that
$$
\nabla  \left( \Rrhom H A {\oaun 1} M\right) \subseteq M^{coH},
$$
 then  $(M^{co H},\nabla^{co H})$ is a hom-connection with respect to $\Omega B$. 
\end{corollary}

\begin{proof}
Let $s:A\rightarrow B\otimes A$ be a normalised left $B$-linear, right $H$-colinear splitting of the multiplication $\mu_B: B\otimes A \rightarrow A$. 
Define the map $\Pi:\oaun 1 \rightarrow (\obun 1)A$ by 
$$
\Pi (d(a)a')= s(a)a' - 1\otimes aa',
$$
 for all $a,a'\in A$.
Since $\mu_B(\Pi(d(a)a'))  = \mu_B(s(a)a') - aa' = 0$ the map is well defined as $\ker(\mu_B)=(\obun 1)A$. Obviously $\Pi$ is right $A$-linear. Since $s$ is left $B$-linear, $\Pi$ is also left $B$-linear, because, for all $a\in A$, $b\in B$, 
$$
\Pi(bd(a))=\Pi(d(ba)-d(b)a) = s(ba)-1\otimes ba - s(b)a + 1\otimes ba = b\Pi(da).$$ 
Moreover $\Pi$ is right $H$-colinear since $s$ is:
\begin{eqnarray*}
(\Pi\otimes \id)\circ \varrho^{\oaun 1}(da) &=& \Pi(da \sw 0)\otimes a \sw 1 = s(a\sw 0)\otimes a\sw 1 - 1\otimes a\sw 0\otimes a\sw 1 \\
&=& \varrho^{\oaun 1}\left(s(a) - 1\otimes a\right)  = \varrho^{\oaun 1}(\Pi(da)).
\end{eqnarray*}
Take an element of the form $d(b)a$ with $b\in B, a\in A$, then 
$$
\Pi(d(b)a) = s(b)a - 1\otimes ba = b(1\otimes 1)a - 1\otimes ba = d(b)a, 
$$
as $s(1)=1\otimes 1$. Hence $\Pi$ splits the inclusion $(\obun 1)A \subseteq \oaun 1$ as a left $B$-, right $A$-module and right $H$-colinear map.

Since a principal comodule algebra is a faithfully flat Hopf-Galois extension, the final assertion follows by Theorem~\ref{thm.induction}.
\end{proof}

The induction procedure of hom-connections presented in Theorem~\ref{thm.induction} can be performed for $A=\cO_q(SL(2))$. As explained in Section~\ref{sec.sl2}, $\cO_q(SL(2))$ is a $\Z$-graded algebra. The degree zero subalgebra is generated by $\alpha\beta$, $\gamma\delta$ and $\beta\gamma$ and is known as the algebra of functions on the {\em standard quantum} or {\em  Podle\'s sphere} $\cO_q(S^2)$ \cite{Pod:sph}. The statement: $\cO_q(SL(2))$ is a $\Z$-graded algebra can be rephrased equivalently as: $\cO_q(SL(2))$ is a comodule algebra over the Hopf algebra $H=\k \Z = \k [z,z^{-1}]$ (the group algebra of $\Z$ or the algebra of Laurent polynomials in one variable). In this set-up the algebra of functions on the quantum sphere is the fixed point (coinavriant) subalgebra of $\cO_q(SL(2))$, i.e.\ $\cO_q(S^2) =\cO_q(SL(2))^{co H}$. Furthermore, $\cO_q(SL(2))$ is a {\em strongly} graded algebra (meaning that $A_kA_l = A_{k+l}$, for degree $k$, $l$, $k+l$ parts of $A$).  In the Hopf-algebraic terms this means that $\cO_q(S^2) \subseteq \cO_q(SL(2))$ is a Hopf-Galois extension by $H = \k [z,z^{-1}]$. In fact, $\cO_q(SL(2))$ is a principal comodule algebra. It is a $q$-deformation of the classical Hopf fibration and one of the first examples of quantum principal bundles; see \cite{BrzMaj:gau}. 

Take $\oan 1$ to be the 3D calculus described in Section~\ref{sec.sl2}. This can be seen to induce the calculus $\obn 1$ on $\cO_q(S^2)$ as follows (see \cite{Maj:Rie} for a detailed description). First note that the $\Z$-grading of $\cO_q(SL(2))$ can be extended to $\oan 1$ by setting $|\omega_-| = -2$, $|\omega_0| = 0$, $|\omega_+| = 2$. Then the exterior differential is a degree-preserving map. This means that there is a coaction of $H$ on $\oan 1$ compatible with the $A$-mutliplication and with $d$, i.e.\ $\oan 1$ is a covariant calculus.  An easy calculation yields  
\begin{eqnarray*}
d (\alpha\beta) &=& \alpha^2\,\omega_-  - q^2\,\beta^2\,\omega_+ = q^2 \omega_-\alpha^2\,  - \omega_+\, \beta^2\, ,\cr
q\,d (\beta\gamma) &=& \alpha\gamma\,\omega_-  - 
q^2\,\beta\delta\,\omega_+ = q^2 \omega_-\alpha\gamma\,  - 
\omega_+\beta\delta \ ,\cr
d (\gamma\delta) &=& \gamma^2\,\omega_-  - q^2\,\delta^2\,\omega_+ = q^2 \omega_-\gamma^2\,  - \omega_+ \delta^2\ .
\end{eqnarray*}
Thus $\obn 1$ is a left and right $B$-module generated by 
$$
\alpha^2\,\omega_-, \quad \alpha\gamma\,\omega_-, \quad \gamma^2\,\omega_-, \quad \beta^2\,\omega_+, \quad \beta\delta\,\omega_+, \quad \delta^2\,\omega_+.
$$
The quantum determinant relations (see equations \eqref{q.det} below)
imply that $\obn 1A$ is a free $A$-module generated by $\omega_-$ and $\omega_+$. Since $\oan 1$ is a free module generated by $\omega_0, \omega_\pm$, there is a right $A$-module map
$$
\bar\Pi : \oan 1 \to \obn 1 A, \qquad \omega_-a_- + \omega_0 a_0 + \omega_+ a_+ \mapsto \omega_-a_- + \omega_+ a_+,
$$
which splits the inclusion $\obn 1A\subset \oan 1$. The map $\bar\Pi$ preserves $\Z$-grades, hence it is a right $H$-comodule map. It is also clearly a left $B$-module map (note that $\alpha\beta$, $\gamma\delta$ and $\beta\gamma$ commute with $\omega_\pm$). 

For the $q$-deformed Hopf bundle, the translation map $\tau: H\to A\ot_B A$, $h\mapsto \can^{-1}(1\ot h)$ is given by $\tau(h) = S(\mathfrak{i}(h)\sw 1) \ot_B \mathfrak{i}(h)\sw 2$, where the map $\mathfrak{i}: H    \to A$ is given as follows:  
$$
\mathfrak{i}(1) = 1, \qquad \mathfrak{i}(z^n) =\alpha^n, \qquad \mathfrak{i}(z^{-n}) = \delta^n, \qquad \mbox{for all } n\in \N .
$$
Define
$
\Pi: \oan 1 \to \obn 1 \ot_B A$ by setting
$$
 \Pi(\omega) = \bar\Pi(\omega S(\mathfrak{i}(z^n)\sw 1))\ot_B\mathfrak{i}(z^n)\sw 2,
$$
for any $\omega \in \oan 1$ of $\Z$-degree $n$. Thus, for all $\omega \in \oan 1$, $\Pi(\omega) = \bar\Pi(\omega\sw 0)\tau(\omega\sw 1)$.  Since $\bar\Pi$ is left $B$-linear, so is $\Pi$. The right $H$-linearity of $\Pi$ follows from the right $H$-linearity of $\bar\Pi$ and $\mathfrak{i}$. To prove the right $A$-linearity of $\Pi$ one uses the right $A$-linearity of $\bar\Pi$ together with the following two properties of the translation map. For all $b\in B$ and $g,h\in H$, $b\tau(h) = \tau(h) b$ and $\tau(gh) = h\su 1\tau(g)h\su 2$, where the notation $\tau(h) = h\su 1\ot_B h\su 2$ (summation implicit) is used; see \cite[3.4]{Sch:rep}. The map $\Pi$ splits the map $j: \obn 1\ot_B A\to \oan 1$.

We are now in position to apply Theorem~\ref{thm.induction} to hom-connection $(\nabla, A)$ given in equation \eqref{eq.hom.3d}. The statement $f\in \Rrhom HA {\oan 1} A$ means that $f$ is a $\Z$-degree preserving map. Thus $|f(\omega_0)| =0$, $|f(\omega_-)| =-2$, $|f(\omega_+)| =2$. The definitions of $q$-skew derivations $\partial_0$, $\partial_\pm$ imply that $\partial_0$ is a degree zero map, $\partial_+$ lowers degree by two, while $\partial_-$ raises degree by two. Therefore, for any $f\in \Rrhom HA {\oan 1} A$, $|\nabla (f)| =0$. This means that $\nabla(f) \in B = \cO_q(S^2)$, i.e.\ $\nabla$ satisfies the requirements of Theorem~\ref{thm.induction}. Thus the hom-connection \eqref{eq.hom.3d} induces a hom-connection $(\nabla^{coH},  \cO_q(S^2))$, for all $f \in \rhom B {\obn 1} B$ given by
\begin{equation}\label{eq.con.s2}
\nabla^{coH}(f) =  q^{-2}\partial_+\left( \widehat{f}\left(\omega_+\right)\right) + q^{2}\partial_-\left( \widehat{f}\left(\omega_-\right)\right),
\end{equation}
where 
$$
\widehat{f}(\omega) = f(\omega S(\mathfrak{i}(z^n)\sw 1))\mathfrak{i}(z^n)\sw 2, 
$$
for any $\omega \in \obn 1A$ of $\Z$-degree $n$.
Noting that the $\omega_\pm$ have $\Z$-degrees $\pm 2$ and taking into account the formulae \eqref{cop.sl2} for the coproduct  on and commutation rules \eqref{comm.sl2} in $\cO_q(SL(2))$, the definition of twisted derivations $\partial_\pm$  \eqref{eq.partial.pm} and corresponding automorphisms $\sigma_\pm$   \eqref{sigma.3d}, a straightforward calculation yields the following explicit form of $\nabla^{coH}$:
\begin{eqnarray}\label{hom.con.b}
\nabla^{coH}(f) \!\!\! &= & \!\!\! q^2\partial_-\left(f\left(\alpha^2 \omega_-\right)\right)\delta^2 - (q^3+ q) \partial_-\left(f\left(\alpha\gamma \omega_-\right)\right)\beta\delta 
+ q^4\partial_-\left(f\left(\gamma^2 \omega_-\right)\right)\beta^2 \nonumber \\
&& \!\!\! \!\!\! \!\!\! \!\!\!+ q^{-4}\partial_+\left(f\left(\beta^2 \omega_+\right)\right)\gamma^2
- (q^{-3}+ q^{-5}) \partial_+\left(f\left(\beta\delta \omega_+\right)\right)\alpha\gamma 
+ q^{-2}\partial_+\left(f\left(\delta^2 \omega_+\right)\right)\alpha^2 \nonumber \\
&& \!\!\! \!\!\! \!\!\! \!\!\! +(q+ q^{-1})\left[\left(q^2f\left(\gamma^2\omega_-\right) - f\left(\delta^2 \omega_+\right)\right)\alpha\beta + \left(f\left(\alpha^2\omega_-\right) - q^{-2}f\left(\beta^2 \omega_+\right)\right)\gamma\delta\right. \nonumber \\
&&\;\;\;\;\;\;\; \left.- \left(qf\left(\alpha\gamma\omega_-\right) - q^{-1}f\left(\beta\delta \omega_+\right)\right)\left(\alpha\delta +q^{-1}\beta\gamma\right)\right]\, .
\end{eqnarray}

We now proceed to identify integral forms associated to $\nabla^{coH}$ with the de Rham complex of $\cO_q(S^2)$.  The higher-order differential calculus relations \eqref{eq.3Drel.higher}, \eqref{eq.3Drel.d.higher}  restrict to produce the higher order differential calculus on $\cO_q(S^2)$; see \cite{Maj:Rie}. The module of two-forms $\obn 2$ is freely generated by the central element  $\omega_+\omega_- =- q^2\omega_-\omega_+$. Any $\omega \in \obn 1$ can be written as $\omega = x\omega_- +y\omega_+$, for some $x,y \in \cO_q(SL(2))$ of $\Z$-degrees $|x| =2$, $|y| = -2$. The differential on such an $\omega$ is given by
\begin{equation}\label{db1}
d\omega = \left(\partial_+(x) - q^{-2}\partial_-(y)\right) \omega_+\omega_-.
\end{equation}
The hom-connection \eqref{hom.con.b} can be identified with the differential $d: \obn 1 \to \obn 2$ as follows. Set 
 $$
 \sfb_1:=\alpha^2,\,\, \sfb_2:=\gamma^2,\,\, \sfb_3:=\alpha\gamma; \quad \sfa_1:=\delta^2, \,\, \sfa_2:=q^2\beta^2,\,\, \sfa_3:= -(q+ q^{-1})\beta \delta \ ,
 $$ and  
 $$ q_1:= 1,\,\, q_2:=q^{-4}, \,\, q_3 := q^{-2} .
 $$
Then the quantum determinant condition yields the following equalities in $B=\cO (S_q^2)$,
 \begin{equation}\label{q.det} 
 \sum_i \sfb_i\sfa_i \,\, =\,\, \sum_i q_i\sfa_i \sfb_i \,\, =\,\, 1.
\end{equation} 
$\obn 1$ decomposes into the direct sum $\obn 1 = \Rp \oplus \Rm$, where $\Rp$ is generated by $\sfw_i:= \omega_+ \sfa_i$ and $\Rm$ is generated by $\sfu_i:= \omega_-\sfb_i$.  Each of the $\Gamma_\pm$  is a finitely generated and projective module with the respective dual bases $\sfw_i^*\in \Rp^* := \rhom B {\Rp} B$, $\sfu_i^*\in \Rm^* := \rhom B {\Rm} B$, given by
 \begin{equation}\label{dual.basis}
 \sfw_i^*(\sfw_j) = q_i\sfb_i\sfa_j, \qquad  \sfu_i^*(\sfu_j) = \sfa_i\sfb_j, \qquad i,j=1,2,3.
\end{equation}
 On these dual basis generators the hom-connection $\nabla^{coH}$ of equation \eqref{hom.con.b} comes out as
 \begin{equation}\label{hom.dual}
  \nabla^{coH}\left( \sfw_i^*\right) = q_iq^{-2}\partial_+(\sfb_i), \qquad 
   \nabla^{coH}\left( \sfu_i^*\right) = q^{2}\partial_-(\sfa_i), \qquad i=1,2,3.
\end{equation}
 Using dual bases one constructs the following $B$-bimodule isomorphisms
 $$
 \xymatrix@R=0pt{ \Psi_{+}: \Rp \ar@{->}[rr] & & \Rm^*, \\ ~\hspace{.2in}\sfw \ar@{|->}[rr] & & \sum_i \sfu_i^* \sfw_i^*(\sfw),  }\qquad 
 \xymatrix@R=0pt{ \Psi_{-}: \Rm \ar@{->}[rr] & & \Rp^*,  \\ ~\hspace{.2in}\sfu \ar@{|->}[rr] & & \sum_i q_i^{-1} \sfw_i^* \sfu_i^*(\sfu) .}
 $$
 These combine into an isomorphism of $B$-bimodules
 $$
 \Psi := \Psi_+ - q^2\Psi_-: \Rp\oplus \Rm = \obn 1 \longrightarrow \obn 1^* = \Rm^*\oplus \Rp^*.
 $$
 Let $\Theta: B\to \obn 2$ be the isomorphism given by $b\mapsto \omega_-\omega_+ b$. A straightforward calculation that uses  definitions of $\Theta$ and $\Psi$ as well as formulae \eqref{q.det}, \eqref{dual.basis} and  \eqref{hom.dual} yields
 $$
 d = \Theta \circ \nabla^{coH}\circ \Psi,
 $$
 where $d$ is the differential described in \eqref{db1}.  
 
Let $\phi: \obn 2\to B$ be the right $B$-linear map dual to $\omega_-\omega_+$. Equations \eqref{q.det},  \eqref{dual.basis}  and \eqref{hom.dual} allow one to compute that $\nabla^{coH}_1(\phi) =0$. In view of the isomorphism $\Theta^* : B \to \obn 2^*$, $b\mapsto \phi b$ and right $B$-linearity of the curvature, the hom-connection $\nabla^{coH}$ is flat. Similarly to Proposition~\ref{prop.3d.flat}, the preceding discussion yields
\begin{proposition}\label{prop.sphere}
Let $B=\cO_q(S^2)$ and $\ob$ be the differential graded algebra described above. The 
diagram 
$$
\xymatrix{ B\ar[rr]^d \ar[d]_{\Theta^*} && \obn 1\ar[rr]^d \ar[d]_\Psi  &&\obn 2\\
\obn 2^* \ar[rr]^{\nabla^{coH}_1}  && \obn 1^*\ar[rr]^{\nabla^{coH}} && B \ar[u]_\Theta \, ,}
$$
in which all columns are ($B$-bimodule) isomorphisms, is commutative. Consequently, the integral complex associated to $\nabla^{coH}$ is isomorphic to the de Rham complex.
\end{proposition}
The $\nabla^{coH}$-integral $\Lambda^{coH}$ on $\cO_q(S^2)$ is then obtained as the restriction of the $\nabla$-integral $\Lambda$ on $\cO_q(SL(2))$ to $\cO_q(S^2)$; see Section~\ref{sec.sl2}. Therefore, up to a scalar multiple, $\Lambda^{coH}$ is equal to the unique normalised $\cO_q(SL(2))$-invariant functional on $\cO_q(S^2)$ described in \cite{NouMim:sph}. The only non-zero values of $\Lambda^{coH}$ are the same as those in \eqref{integral}.

\setcounter{equation}{0}

  \section*{Acknowledgements} 
T.\ Brzezi\'nski would like to thank Theodore Voronov for drawing his attention to reference \cite{Man:gau}.  L.\ El Kaoutit and C.\ Lomp would like to thank all the members of Mathematics Department at Swansea University for a warm hospitality. The stay of L.\ El Kaoutit was supported by the grant JC-2008-00440 from the Ministerio De Ciencia E Innovaci\'on of Spain.
The stay of C.\ Lomp was supported by the grant SFRH/BSAB/826/2008 from Funda\c{c}\~ao para a Ci\^encia e a Tecnologia (Portugal).

\end{document}